\documentclass[MSNbibl,number,citesort,seceqn,dvips]{arxbj}
\usepackage{upgreek}

\aid{0}
\volume{20}
\issue{2}
\pubyear{2014}
\firstpage{747}
\lastpage{774}
\doi{10.3150/12-BEJ505} 

\makeatletter

\newtheorem{lemmma}{Lemma}[section]

\newcommand{\RMo}{\mathrm{o}}
\newcommand{\RMO}{\mathrm{O}}

\newcommand{\RMe}{\mathrm{e}}

\newcommand{\mrmd}{\,\mathrm{d}}
\newcommand{\dd}{\mathrm{d}}

\newtheorem{theorem}{Theorem}
\newtheorem{lemma}{Lemma}
\newremark{remark}{Remark}
\newtheorem{corollary}{Corollary}
\newproclaim{assumption}{Assumption}

\newcommand{\dct}{\operatorname{DCT}_8}

\newcommand{\R}{\mathbb{R}}
\newcommand{\N}{\mathbb{N}}

\newcommand{\Var}{\operatorname{Var}}
\newcommand{\Cov}{\operatorname{Cov}}
\newcommand{\E}{\mathbb{E}}
\renewcommand{\P}{\mathbb{P}}

\newcommand{\MSE}{\operatorname{MSE}}

\newcommand{\signtxt}{\operatorname{sign}}
\renewcommand{\limsup}{\mathop{\overline{\lim}}}
\renewcommand{\liminf}{\mathop{\underline{\lim}}}
\newcommand{\oracle}{\mathrm{oracle}}
\newcommand{\id}{\operatorname{id}}
\newcommand{\GM}{\operatorname{GM}}
\newcommand{\sigmamin}{\underline{\sigma}}
\newcommand{\sigmamax}{\overline{\sigma}}

\newcommand{\bR}{\mathbf{R}}
\newcommand{\bU}{\mathbf{U}}
\newcommand{\bX}{\mathbf{X}}
\newcommand{\bY}{\mathbf{Y}}
\newcommand{\bZ}{\mathbf{Z}}

\makeatother

\begin{document}
\begin{frontmatter}

\title{Asymptotically efficient estimation of a scale parameter in
Gaussian time series and closed-form expressions for the Fisher information}
\runtitle{Estimation of a scale parameter}

\begin{aug}
\author[1]{\fnms{Till} \snm{Sabel}\thanksref{1}\ead[label=e1]{tsabel@uni-goettingen.de}} \and
\author[2]{\fnms{Johannes} \snm{Schmidt-Hieber}\corref{}\thanksref{2}\ead[label=e2]{Johannes.Schmidt.Hieber@ensae.fr}}
\runauthor{T. Sabel and J. Schmidt-Hieber} 
\address[1]{Institut f\"ur Mathematische Stochastik,
Universit\"at G\"ottingen, Goldschmidtstr. 7, 37077 G\"ottingen,
Germany. \printead{e1}}
\address[2]{CREST-ENSAE, 3, Avenue Pierre Larousse, 92240 Malakoff,
France.\\ \printead{e2}}
\end{aug}

\received{\smonth{8} \syear{2012}}
\revised{\smonth{11} \syear{2012}}

%
\begin{abstract}
Mimicking the maximum likelihood estimator, we construct first order
Cramer--Rao efficient and explicitly computable estimators for the scale
parameter $\sigma^2$ in the model $Z_{i,n}=\sigma n^{-\beta}X_{i}+
Y_{i}, i=1,\ldots,n, \beta>0$ with independent, stationary Gaussian
processes $(X_i)_{i\in\mathbb{N}}$, $(Y_i)_{i\in\mathbb{N}}$, and
$(X_i)_{i\in\mathbb{N}}$ exhibits possibly long-range dependence. In a
second part, closed-form expressions for the asymptotic behavior of the
corresponding Fisher information are derived. Our main finding is that
depending on the behavior of the spectral densities at zero, the Fisher
information has asymptotically two different scaling regimes, which are
separated by a sharp phase transition. The most prominent example
included in our analysis is the Fisher information for the scaling
factor of a high-frequency sample of fractional Brownian motion under
additive noise.
\end{abstract}

%

\end{frontmatter}

\section{Introduction}\label{intro}

Let $(X_i)_{i\in\mathbb{N}}$ and $(Y_i)_{i\in\mathbb{N}}$ be
independent Gaussian processes with known distribution. Suppose that
we observe $\bZ:=\bZ_n:=(Z_{1,n},\ldots,Z_{n,n})$ with
%
\begin{equation}\label{eqmod}
Z_{i,n}=\sigma n^{-\beta} X_i+ Y_i,\qquad
i=1,\ldots,n\mbox{ and } \beta, \sigma>0.
\end{equation}
In our framework, the parameter $\beta$ is assumed to be known. We are
interested in the case where $(X_i)_{i\in
\mathbb{N}}$ is stationary and $(Y_i)_{i\in\mathbb{N}}$ is a noise
process. Our theory includes white noise and
increments of white noise as special cases for $(Y_i)_{i\in\mathbb
{N}}$ (cf. Assumption \ref{assY}). The problem, which we address in
this work, is asymptotically optimal
estimation of the scale parameter $\sigma^2$. In order to understand
its asymptotic properties, the key ingredient is knowledge of the
Fisher information, for which closed-form expressions will be derived
as well.

Our study is motivated by estimation of the variance $\sigma^2$ of a
fractional Brownian
motion (fBM) $(B_t^H)_{t\geq0}$ at time points $i/n, i=1,\ldots,n$,
under additive Gaussian white
noise (WN), that is,
%
\begin{equation}\label{eqmodfBM}
V_{i,n}:=\sigma B^H_{i/n}+\tau
\varepsilon_{i},\qquad i=1,\ldots,n.
\end{equation}
Here, $H$ refers to the Hurst index (or self-similarity parameter) and
$(\varepsilon_{i})_{i}$ is a sequence of i.i.d. standard normal random variables.
This model has attracted a lot of attention, recently (cf. Gloter and
Hoffmann \cite{GloterHoffmann2004,GloterHoffmann2007} and for the special case $H=1/2$, cf. Stein \cite
{ste}, Gloter and Jacod \cite{GloterJacod2001a,GloterJacod2001b},
as well as Cai \textit{et al.}  \cite{CaiMunkSchmidtHieber2010}). Let us
call it the fBM${}+{}$WN model and note that the
increment vector is of type (\ref{eqmod}) with $\beta=H$. This shows
that models (\ref{eqmod}) and
(\ref{eqmodfBM}) coincide, if $X_i$ and $Y_i$ are chosen as the
increments $n^H(B^H_{i/n}-B^H_{(i-1)/n})$ and
$\tau(\varepsilon_{i,n}-\varepsilon_{i-1,n})$, with $\varepsilon_{0,n}:=0$,
respectively.

Estimation of $\sigma^2$ (and $H$) was discussed in slightly more
general settings than the fBM${}+{}$WN model by Gloter and Hoffmann
\cite{GloterHoffmann2004,GloterHoffmann2007}. In these papers, it was
proven that for $H>\frac{1}{2}$ the optimal rate of convergence for
$\sigma^2$ is $n^{-1/(4H+2)}$. More extensively studied and of
particular interest is the case $H=\frac{1}{2}$, due to its
applications to high-frequency modeling of stock returns. For this
case, the asymptotic Fisher information is known to be
$n^{1/2}(8\tau\sigma^3)^{-1} (1+\RMo(1) )$ (cf. Gloter and Jacod
\cite{GloterJacod2001a,GloterJacod2001b}, and Cai \textit{et al.}
\cite {CaiMunkSchmidtHieber2010}). This result had a big impact as a
benchmark for estimation of the integrated volatility (cf.
Barndorff-Nielsen \textit{et al.}
\cite{Barndorff-NielsenHansenLundeShephard2008}, Podolskij and Vetter
\cite{PodolskijVetter2009}, Jacod \textit{et
al.}~\cite{JacodLiMyklandPodolskijVetter2009}, and Zhang \cite
{Zhang2006}) as well as for the asymptotic equivalence theorem by
Rei\ss\ \cite{Reiss2011}. The fact that the multiplicative inverse of
the asymptotic Fisher information is linear in $\tau$ and proportional
to the cube of $\sigma$ is surprising and requires further
understanding.

The main contribution of our work to the existing literature is that
for $0<H<1$, the Fisher information $I_{\sigma^2}^n$ for estimation of
$\sigma^2$ in the fBM${}+{}$WN model is given by
%
\begin{equation}\label{eqFishInfo}
I_{\sigma^2}^n= n^{1/({2H+1})} \sigma^{-({8H+2})/({2H+1})}
\tau^{-2/({2H+1})}c_H+\RMo \bigl(n^{1/({2H+1})} \bigr),
\end{equation}
where $c_H$ is a constant only depending on $H$ (for an explicit
expression of $c_H$, cf. Corollary~\ref{cortheoremholdsforfBM}).

In general, we focus on the situation, where the Fisher information
converges to infinity for $n\rightarrow\infty$, which corresponds to
consistent estimation of $\sigma^2$. In view of $n^{-\beta}
X_{i}=\RMO_p(n^{-\beta})$ and $Y_i=\RMO_p(1)$, it is not clear at all that
there are such situations. In fact, the rate at which the Fisher
information tends to infinity can be rather unexpected. In a first
place, one might guess that the optimal rate of convergence for
estimation of the ``parameter'' $\sigma^2n^{-2\beta}$ is $n^{-
{1/2}}$ and
hence the Fisher information of $\sigma^2$ should be of the order
$n^{1-4\beta}$ (corresponding to the rate of convergence $n^{2\beta
-1/2}$). However, this heuristic reasoning is in general not true and
\textit{better rates} can be obtained, as, for instance, in (\ref
{eqFishInfo}). Surprisingly, the asymptotic Fisher information has two
different scaling regimes. In fact we will see that for any pair
$(X_i)_i$ and $(Y_i)_i$ there is a positive characteristic $\Diamond$
such that (up to sub-polynomial factors) $I_{\sigma^2}^n\propto
n^{1-\Diamond\beta}$ if $\Diamond<4$ and $I_{\sigma^2}^n\propto
n^{1-4\beta}$ if $\Diamond\geq4$. The latter appears to be the same
rate as in our heuristic argument above. Altogether, the different
scaling behavior becomes visible as elbow effect in the convergence
rate of $\sigma^2$. As a curious fact, let us mention that the
spectral densities of the processes do not need to be known explicitly
in order to compute the proposed estimator or the asymptotic Fisher information.

It is a classical result that if we observe a sample of a stationary
Gaussian process with a spectral density $h(\theta,\cdot)$, the
asymptotic Fisher information $I_{\theta}^n$ for estimation of a
one-dimensional parameter $\theta$ is given by (cf. Davies \cite
{Davies1973} and Dzhaparidze \cite{Dzhaparidze1986} for the general
case as well as Fox and Taqqu \cite{FoxTaqqu1986}, Dahlhaus \cite
{Dahlhaus1989}, Giraitis and Surgailis \cite{GiraitisSurgailis1990}
for long-range dependent processes)
%
\begin{equation}\label{eqFInfogen}
I_\theta^n = \frac{n}{2\uppi } \int_0^\uppi
\bigl(\partial_\theta\log h(\theta,\lambda) \bigr)^2 \mrmd
\lambda+\RMo(n).
\end{equation}
In Theorem \ref{mainthm}, we prove that under fairly general
conditions on $(X_i)_i$, a result of the type (\ref{eqFInfogen})
holds for $\theta=\sigma^2$ in model (\ref{eqmod}). One should
notice that our setting is nonstandard and not covered within the
existing literature. In contrast, due to the factor $n^{-\beta}$, we
cannot work with a fixed $h$ but rather have to consider a sequence of
spectral densities $(h_n)_n$ with degenerate limit. Furthermore, we are
not in the classical parametric estimation setting, that is, $I_{\sigma
^2}^n$ may diverge with a rate which is much slower than $n$. As, for
example, in (\ref{eqFishInfo}), we need therefore to prove (\ref
{eqFInfogen}) with an approximation error which is of smaller order
than $\RMo(n)$. This in turn implies that very precise control on the
(large) noise process $(Y_i)_i$ has to be imposed (cf. Assumption \ref
{assY}). Let us also mention that we cover both cases, long and
short-range dependence of $(X_i)_i$. In particular, this allows to
treat model (\ref{eqmodfBM}) for all $H\in(0,1)$.

The work is organized as follows. In Section \ref{subsecconstruct},
we construct the estimator and investigate its theoretical properties.
Closed-form expressions for the Fisher information are derived in
Section~\ref{secclosedform}. In particular, we give some heuristic
arguments why different scaling regimes appear. To illustrate the
results, some examples are provided in Section \ref{secexamples}.
Proofs are deferred to the \hyperref[app]{Appendix} and the Supplement \cite
{SabelSchmidtHieber}.

\textit{Notation}: We write $\bX:=\bX_n:=(X_1,\ldots,X_n)$, $\bY:=\bY
_n:=(Y_1,\ldots,Y_n)$ and $\bZ:=\bZ_n:=(Z_{1,n},\ldots,Z_{n,n})$.
For two sequence $(a_k)_k$ and $(b_k)_k$, we say that $a_k\sim b_k$ iff
$\lim_{k\rightarrow\infty} a_k/b_k= 1$. Similar, for two functions
$g_1$ and $g_2$, we write $g_1(\lambda)\sim g_2(\lambda)$ (for
$\lambda\downarrow0$) iff $\lim_{\lambda\downarrow0} g_1(\lambda
)/\allowbreak g_2(\lambda) =1$.

\section{Main results}

Let $\bU=\bU_n$ be an $n$-dimensional, centered Gaussian vector with
positive definite covariance matrix $\Sigma_\theta$, depending on a
one-dimensional parameter $\theta\in\mathbb{R}$. The log-likelihood
function is
\[
L(u|\theta) = - \frac n2 \log(2\uppi ) -\frac12 \log\bigl(|\Sigma_\theta|\bigr) -
\frac12 u^t \Sigma_\theta^{-1} u
\]
with $|\Sigma_\theta|$ the determinant of $\Sigma_\theta$. Let
$\partial_\theta\Sigma_\theta$ denote the
entrywise derivative of $\Sigma_\theta$ with respect to $\theta$
(which we assume to exist). Since $\partial_\theta\log(|\Sigma
_\theta|)=
\operatorname{tr}(\Sigma_\theta^{-1}\partial_\theta\Sigma_\theta
)$ and
$\partial_\theta\Sigma_\theta^{-1}= - \Sigma_\theta^{-1}
(\partial_\theta\Sigma_\theta) \Sigma_\theta^{-1}$, we find for
the score function
\[
\dot{L}(\bU|\theta):=\partial_\theta L(\bU|\theta) = -\tfrac12
\operatorname{tr}\bigl(\Sigma_\theta^{-1}\partial_\theta
\Sigma_\theta\bigr) + \tfrac12 \bU^t \Sigma_\theta^{-1}
(\partial_\theta\Sigma_\theta) \Sigma_\theta^{-1}
\bU.
\]
In distribution, $\bU=\Sigma_\theta^{1/2} \xi$ for an
$n$-dimensional standard normal vector $\xi$. Together with some
algebra this shows that the Fisher information for $\theta$ is
$I_\theta^n=\frac12 \operatorname{tr} ([(\partial_\theta
\Sigma_\theta)\Sigma_\theta^{-1}]^2 )$ (cf. also Porat and
Friedlander \cite{PoratFriedlander1986}). In particular, for model
(\ref{eqmod}) we obtain
%
\begin{equation}\label{eqFishexplexpr}
I_{\sigma^2}^n = \tfrac12 \operatorname{tr} \bigl(
\bigl[n^{-2\beta}\Cov(\bX)\Cov(\bZ)^{-1}\bigr]^2
\bigr).
\end{equation}
To simplify the notation, we will view $I_{\sigma^2}^n$ in the
following always as a sequence in $n$.

\subsection{An asymptotically Cramer--Rao efficient estimator}
\label{subsecconstruct}

In this section, we construct an explicitly computable estimator which
mimics the MLE. Furthermore, we prove that the mean squared error (MSE)
of this estimator is first order optimal (cf. Theorem \ref
{thmestisopt}). All the results in this section work under fairly
general conditions. In fact, we only require that $\Cov(\bX)$ and
$\Cov(\bY)$ are known and positive definite for all $n$ as well as
divergence of the Fisher information. In particular, we neither have to
impose stationarity on $(X_i)_i$ or $(Y_i)_i$ nor do we assume that
$\Cov(\bX)$ and $\Cov(\bY)$ have the same set of eigenvectors.

The construction will be done in several steps. First, we can find
$n\times n$ matrices $A$ and $D$ such that $\Cov(\bY)=A^tA$,
$D^tD=\id_n$ (the identity), and $D$ diagonalizes $(A^{-1})^t \Cov
(\bX) A^{-1}$. Hence, $\Lambda:= (A^{-1}D)^t
\Cov(\bX) (A^{-1}D)$ is diagonal and the diagonal entries are denoted
by $\lambda_1,\ldots,\lambda_n$. The
maximum likelihood equation in the transformed model $(\widetilde
Z_{1,n},\ldots, \widetilde Z_{n,n})^t:=
(A^{-1}D)^t\bZ$ motivates to consider the oracle estimator
%
\begin{equation}\label{eqoracleestsigmasqr}
\widehat\sigma_{\oracle}^2:= \bigl(2I_{\sigma^2}^n
\bigr)^{-1} \sum_{i=1}^n
\frac{\lambda_i
n^{-2\beta}(\widetilde
Z_{i,n}^2-1)}{(\sigma^2n^{-2\beta}\lambda_i+1)^2} = \sigma^2+\bigl
(2I_{\sigma^2}^n
\bigr)^{-1}\sum_{i=1}^n
\frac{\lambda_i n^{-2\beta} (\widetilde Z_{i,n}^2-\E\widetilde
Z_{i,n}^2)}{(\sigma^2n^{-2\beta} \lambda_i+1)^2}.\quad
\end{equation}
To verify the equality, one should note that by rewriting (\ref
{eqFishexplexpr})
\[
I_{\sigma^2}^n = \frac12 \sum_{i=1}^n
\frac{\lambda_i^2n^{-4\beta
}}{(\sigma^2\lambda_in^{-2\beta}+1)^2}.
\]
Observe that the oracle estimator (\ref{eqoracleestsigmasqr}) is
unbiased and attains the Cramer--Rao bound since
$\Var(\widehat\sigma_{\oracle}^2)=(I_{\sigma^2}^n)^{-1}$.
Oracle\vspace*{1pt}
estimators depend on the unknown quantities itself,
and are thus not computable. Below, we derive a simple construction for
a statistical estimator which mimics $\widehat\sigma_{\oracle}^2$
and is asymptotically sharp. Similar as in \cite
{CaiMunkSchmidtHieber2010}, we use a sample splitting technique.
First, we take a small part of the data in the transformed model, which
are used for a preliminary estimate, say $\widetilde\sigma^2$, of
$\sigma^2$. In a second step, we plug $\widetilde\sigma^2$ into
(\ref{eqoracleestsigmasqr}). Discarding all indices, which were
already used for $\widetilde\sigma^ 2$ gives an estimator, which as
we show has asymptotically the same properties as $\widehat\sigma
_{\oracle}^2$. This implies then the first order Cramer--Rao efficiency.

The next lemma ensures that sample splitting can be done.
%
\begin{lemma}
\label{lempartitioneigvals}
For $u>0$ and $B\subset\{1,\ldots,n\}$, let
\[
I_u^B:= \frac12 \sum_{i\in B}
\frac{\lambda_i^2n^{-4\beta
}}{(u\lambda_in^{-2\beta}+1)^2}.
\]
Then there is a sequence of index sets $(A_n)_n$ with $A_n \subset\{
1,\ldots,n\}$ such that
\[
I_1^{A_n} \rightarrow\infty\quad\mbox{and}\quad
\frac{I_1^{A_n}}{I_1^n} \rightarrow0 \qquad\mbox{as } n\rightarrow\infty.
\]
\end{lemma}
\begin{pf}
Note that $0\leq c_{i,n}:=\lambda_i^2n^{-4\beta}(n^{-2\beta}\lambda
_i+1)^{-2}<1$. Consider the partial sums $S_k=\sum^k_{i=1}c_{i,n}$
and observe that $S_{k+1}-S_k\leq1$ as well as $S_n=I_1^n\rightarrow
\infty$. Therefore, we can find $k^*(n)$, s.t. $\sqrt{I_1^n}\leq
S_{k^*(n)}\leq\sqrt{I_1^n}+1$. The result follows with $A_n=\{
1,\ldots,k^*(n)\}$.
\end{pf}

Throughout this section, let $(A_n)_n$ be as in the previous lemma and
pick a sequence $(\delta_n)_n$, satisfying $\delta_n\leq1, \delta
_n \rightarrow0$, and $\delta_n^4 I_1^{A_n}\rightarrow\infty$. A
possible choice is $\delta_n = (I_1^{A_n})^{-1/8}$. Observe that
\[
V:= \bigl(2I_1^{A_n} \bigr)^{-1} \sum
_{i\in A_n} \frac{\lambda_i
n^{-2\beta}(\widetilde
Z_{i,n}^2-1)}{(n^{-2\beta}\lambda_i+1)^2} = \sigma^2 +
\bigl(2I_1^{A_n} \bigr)^{-1} \sum
_{i\in A_n} \frac{\lambda_i
n^{-2\beta}(\widetilde
Z_{i,n}^2-\E\widetilde
Z_{i,n}^2)}{(n^{-2\beta}\lambda_i+1)^2}
\]
has expectation $\sigma^2$ and variance bounded by $(\sigma^4\vee1)
(I_{1}^{A_n})^{-1}$. Now, we define the preliminary estimator
$\widetilde\sigma^2$, as the truncated version of $V$,
%
\begin{equation}\label{eqprelimestdef}
\widetilde\sigma^2:= (V \vee\delta_n) \wedge
\delta_n^{-1}.
\end{equation}
This allows us to construct the final estimator $\widehat\sigma_n^2$
for $\sigma^2$. Let $A_n^c= \{1,\ldots,n\}\setminus A_n$ and set
%
\begin{equation}\label{eqestsigmasqr}
\widehat\sigma_n^2:= \bigl(2I^{A_n^c}_{{\widetilde\sigma}^2}
\bigr)^{-1} \sum_{i\in
A_n^{c}} \frac{\lambda_i
n^{-2\beta}(\widetilde
Z_{i,n}^2-1)}{(\widetilde\sigma^2n^{-2\beta}\lambda_i+1)^2}
= \sigma^2+\bigl(2I^{A_n^c}_{\widetilde\sigma^2}
\bigr)^{-1}\sum_{i\in A_n^c} \frac{\lambda_i n^{-2\beta} (\widetilde
Z_{i,n}^2-\E\widetilde
Z_{i,n}^2)}{(\widetilde\sigma^2n^{-2\beta}
\lambda_i+1)^2}.
\end{equation}
One should note the similarity to the oracle $\widehat\sigma_{\oracle
}^2$ as introduced in (\ref{eqoracleestsigmasqr}). As the
following theorem shows, $\widehat\sigma_n^2$ has in fact the same
asymptotic MSE as the oracle, implying Cramer--Rao efficiency.

\begin{theorem}
\label{thmestisopt}
Suppose that the Fisher information diverges and $\Cov(\bY)$ is
positive definite. The estimator $\widehat\sigma_n^2$ defined in
(\ref{eqestsigmasqr}) attains the Cramer--Rao bound asymptotically
over every compact set, not containing zero, that is, for $0<\sigma_{\min
}<\sigma_{\max}<\infty$,
\[
\lim_{n\rightarrow\infty} \sup_{\sigma\in[\sigma_{\min
},\sigma_{\max}]}
I_{\sigma^2}^n\cdot\MSE\bigl(\widehat\sigma_n^2
\bigr)=1.
\]
\end{theorem}

\subsection{Closed-form expressions for the Fisher information}
\label{secclosedform}

So far we have seen that there are estimators which are asymptotically
Cramer--Rao efficient. However, in order to get some understanding of
the asymptotics, we need to study the behavior of the Fisher
information. In this section, we derive explicit closed-form expressions.

To state the results, some definitions, in particular from regular
variation theory, are unavoidable. For the notion of quasi-monotone and
slowly varying functions see the monograph \cite
{BinghamGoldieTeugels1987}. A~positive sequence $(r_j)_j$ is called
O-regularly varying if for any $\lambda> 1$,
\[
0< \liminf_{n \rightarrow\infty} \frac{r_{\lfloor\lambda n\rfloor
}}{r_n} < \limsup
_{n \rightarrow\infty} \frac{r_{\lfloor\lambda n\rfloor
}}{r_n}<\infty
\]
with $\lfloor\cdot\rfloor$ the Gauss bracket. For real sequences
$(a_j)_j$, O-regularly varying quasi-monotonicity is
equivalent to the existence of a positive, nondecreasing, and
O-regularly varying sequence $(r_j)_j$ such that the
sequence $(a_j/r_j)_j$ is decreasing. We say that a sequence $(a_j)_j$
is general monotone if there are
finite constants $C,J_0$, such that for any positive integer $J\geq
J_0$, $\sum_{j=J}^{2J-1} |a_{j+1}-a_j|\leq C |a_J|$. The class of
general monotone sequences will be denoted by $\GM$. It was introduced
and studied recently by Belov \cite{Belov2002}
and Tikhonov \cite{Tikhonov2007}. To simplify some arguments, we have
relaxed the original definition slightly by introducing $J_0$ (this
does not cause any trouble and all results on $\GM$ sequences can be
transferred with obvious changes). The class $\GM$ is fairly general
in the sense that it includes all well-known
generalizations of monotone sequences, such as quasi-monotonicity,
regularly quasi-monotonicity,
O-regularly-quasi-monotonicity and sequences of rest bounded variation.

In order to deal with boundary problems (cf. the second example in
Section \ref{secexamples}), we assume that $(X_i)_i$ is only
approximately stationary in the following sense.

\begin{assumption}[(Assumptions on $\bolds{X}$)]\label{assX}
Suppose that there is a stationary process $(X_i')_i$ and a process
$(R_i)_i$ such that in distribution
\[
X_i =X_i'+R_i\qquad \mbox{for
all } i\in\mathbb{N}
\]
and $(X_i')_i$ and $(R_i)_i$ have the following properties.
\begin{itemize}[(ii)]
\item[(i)] There is a positive and quasi-monotone slowly varying
function $\ell\dvtx (0,\infty)\rightarrow\mathbb{R}^+$, satisfying
%
\begin{equation}
\label{assell} \frac{\ell(x\ell^\kappa(x))}{\ell(x)}
\rightarrow1,\qquad
x\rightarrow\infty \mbox{ for all } \kappa\in\R
\end{equation}
such that for an index $\alpha\in(-1/2,1/2)$,
%
\begin{equation}\label{eqautocorrasympbehavior}
\gamma_k:=\Cov\bigl(X_1',X_{1+k}'
\bigr) \sim\signtxt(-\alpha) k^{-2\alpha
-1} \ell(k) \qquad\mbox{as } k\rightarrow
\infty.
\end{equation}
\item[(ii)] With $\bX'_n:=(X_1',\ldots,X_n')$,
$\bR_n:=(R_1,\ldots,R_n)$ and $\|\cdot\|_2$ the Frobenius norm, we have
the uniform bound
%
\begin{equation}\label{eqRfrobbdcond}
\sup_n \bigl\|\Cov\bigl(\bX_n',
\bR_n\bigr)\bigr\|_2+\bigl\|\Cov(\bR_n)\bigr\|_2
<\infty.
\end{equation}
\end{itemize}
\end{assumption}

Throughout the following, we interpret the autocovariance $(\gamma
_k)_{k\in\mathbb{Z}}$ as a sequence on $\mathbb{Z}$ via $\gamma
_k=\gamma_{-k}$. An example for a (quasi-)monotone slowly varying
$\ell$ is the logarithm $\log(1+\cdot)$. However, (\ref{assell})
does not hold for every slowly varying function. A stronger condition,
which implies (\ref{assell}) is
\[
\lim_{\lambda\rightarrow0} \biggl(\frac{\ell(a \lambda)}{\ell
(\lambda)}-1 \biggr)\log(\lambda)
=0 \qquad\mbox{for an } a>1
\]
(cf. Theorem 1 in \cite{BojanicSeneta1971}). As a consequence (\ref
{assell}) holds whenever $\lim_{\lambda\rightarrow\infty} \ell
(\lambda)\in(0,\infty)$. Moreover, if it is true for $\ell$ then
also for $\ell^\mu$, $\mu\geq0$.

The $n\times n$ matrix $\Delta$ denotes the backward difference
operator, that is,
%
\begin{equation}\label{backwarddiff}
\Delta=\pmatrix{1&
\cr
-1&1
\cr
&\ddots&\ddots
\cr
&&-1&1}.
\end{equation}
For a vector $v=(v_1,\ldots,v_n)^t$, $\Delta
v=(v_1-v_0,v_2-v_1,\ldots,v_n-v_{n-1})^t$ with $v_0:=0$ is the
backwards difference process. Furthermore, the transposed matrix
$\Delta^t$ is the negative forward difference operator $\Delta^t
v=-(v_2-v_1,v_3-v_2,\ldots,v_{n+1}-v_{n})^t$ with $v_{n+1}:=0$. We
assume that for a nonnegative integer $K$, the process $\bY$ is
generated by taking the $K$th finite difference of a white noise
process (alternating between forward and backward differences).

\begin{assumption}\label{assY}
Given a nonnegative integer $K$ and $\tau>0$, assume that $\bY$ is an
$n$-dimensional,
centered Gaussian random vector with covariance matrix
$\tau^2 (\Delta\Delta^t )^{K}$ or $\tau^2 (\Delta
^t\Delta)^{K}$.
\end{assumption}

Assumption \ref{assY} imposes in fact a very serious restriction, but
seems to be somehow unavoidable in order to prove the statement (cf.
also the discussion in the \hyperref[intro]{Introduction}). Our results could be worked
out under more general boundary conditions of the difference operator,
of course. It is indeed
sufficient that $\Cov(\bY)$ can be perfectly diagonalized by a
discrete sine or cosine transform. However, since
the assumption above is somehow the most natural one and allows to
treat the fBM${}+{}$WN model, we will restrict ourselves
to it for sake of simplicity.

Let throughout the paper $f=\sum_{k=-\infty}^\infty\gamma_k \cos
(k\cdot)$ denote the spectral density of $(X_i')_{i\in\N}$. Although
$\bX$ and $\bY$ are stationary only
up to boundary values, we will refer occasionally to
%
\begin{equation}\label{eqsdnotation}
f, 4^K\tau^2\sin^{2K} \biggl(
\frac{\cdot}2 \biggr)\quad \mbox{and}\quad h_n=
\sigma^2n^{-2\beta} f+4^K\tau^2
\sin^{2K} \biggl(\frac
{\cdot}2 \biggr)
\end{equation}
as the spectral density of the processes $\bX$, $\bY$, and $\bZ$,
respectively. Because of the imposed independence of $(X_i)_i$ and
$(Y_i)_i$, $h_n$ is the sum of the spectral densities of $\bX$ and
$\bY$.

Define
%
\begin{equation}
\label{eqdefrn} r_n:= \cases{ n^{1-{\beta}/({K-\alpha})} \bigl(\ell
\bigl(n^{{\beta}/({K-\alpha
})}\bigr) \bigr)^{
1/({2K-2\alpha})}+n^{1-4\beta}, &\quad if $K-
\alpha\neq\frac14$,
\vspace*{2pt}\cr
n^{1-4\beta}\log(n)\ell^2(n), &\quad if $K-
\alpha= \frac14$.}
\end{equation}
Now, we are ready to state the main results of the paper. Surprisingly,
it turns out that the rates depend on $K$ and $\alpha$ only through
their (inverted) difference, that is, the problem is characterized by
\[
\Diamond:=\frac1{K-\alpha}.
\]

\begin{theorem}\label{mainthm}
Work in model (\ref{eqmod}) under Assumptions \ref{assX} and \ref
{assY} and suppose that $K-\alpha>\max\{\beta,(4\alpha+1)\beta,1/4\}$.
Further assume that either:
\begin{enumerate}
\item\label{version1}
$\alpha\in(0,1/2)$, $(\gamma_k)_k$ is O-regularly varying
quasi-monotone,
and $\sum_{k=-\infty}^\infty\gamma_k=0$, or
\item\label{version2}
$\alpha\in(-1/4,0)$ and $(\gamma_k)_k\in\GM$, or
\item\label{version3}
$\alpha\in(-1/2,-1/4]$ and there exists a constant $C_1$, such
that for any $p\in\mathbb{N}$, $|\gamma_{p+1}-\gamma_p|\leq C_1
|\gamma_p|p^{-1}$.
\end{enumerate}
Then, the Fisher information of $\sigma^2$ based on $n$ observations
is given by
%
\begin{equation}\label{fisherintegral}
I_{\sigma^2}^n = \frac{n^{1-4\beta}}{2\uppi }\int^\uppi _0
\frac
{f^2(\lambda)}{h_n^2(\lambda)}\mrmd \lambda\bigl(1+\RMo(1)\bigr)+\RMo(r_n).
\end{equation}
If the condition $K-\alpha>\max\{\beta,(4\alpha+1)\beta,1/4\}$ is
replaced by the weaker assumption $K-\alpha>\max\{\beta,(4\alpha
+1)\beta\}$, imposing additionally $\log(n)\ell^2(n)\rightarrow
\infty$ in the critical case $K-\alpha=1/4$, then (\ref
{fisherintegral}) holds as well, provided there exists a constant
$c_f$, such that
%
\begin{equation}\label{eqaddassumpifDiamondlarge}
\bigl|f(\lambda)-f(\mu) \bigr|\leq c_f \lambda^{2\alpha-2}|\lambda-\mu|
\qquad\mbox{for all } 0<\lambda\leq\mu\leq\uppi .
\end{equation}
\end{theorem}

Let
\[
C(\Diamond,\alpha):= \frac{(2-\Diamond)\Diamond}{8\sin(\Diamond{\uppi }
/{2})} \bigl(2\signtxt(-\alpha) \Gamma(-2
\alpha)\cos(\uppi \alpha) \bigr)^{\Diamond/{2}}.
\]

\begin{theorem}\label{thmexplicit}
Under the assumptions of Theorem \ref{mainthm}, the asymptotic Fisher
information is explicitly given by
%
\begin{equation}\label{eqFisherInfoK-alphalarge}
I_{\sigma^2}^n \sim n^{1-\Diamond\beta} \bigl(\ell
\bigl(n^{\Diamond\beta}\bigr) \bigr)^{\Diamond/2}\sigma^{\Diamond-4}
\tau^{-\Diamond} C(\Diamond, \alpha)\qquad \mbox{if } \Diamond<4
\end{equation}
and
%
\begin{equation}\label{eqFisherInfoK-alphasmall}
I_{\sigma^2}^n \sim\frac{n^{1-4\beta}}{2\uppi \tau^4} \int_0^\uppi
f^2(\lambda) \mrmd \lambda= \frac{n^{1-4\beta}}{2\tau^4}\sum
_{k=-\infty}^\infty\gamma_k^2
\qquad\mbox{if } \Diamond>4.
\end{equation}
For $\ell=|{\log^\rho(\cdot)}|, \rho>-1/2$,
%
\begin{equation}\label{eqFisherInfoK-alphamedium}
I_{\sigma^2}^n \sim n^{1-4\beta}(\log n)^{2\rho+1}
\tau^{-4} \frac{(4\beta)^{2\rho+1}}{2\rho+1}\qquad \mbox{if } \Diamond=4.
\end{equation}
\end{theorem}

\begin{corollary}
\label{cortheoremholdsforfBM}
In the fBM${}+{}$WN model (\ref{eqmodfBM}) with $H\in(0,1)$, it holds that
\[
\Diamond=\frac2{2H+1}<4
\]
and the asymptotic Fisher information for $\sigma^2$ is given by
(\ref{eqFishInfo}) with
\[
c_H:= \frac{H\sin^{1/({2H+1})} (\uppi  H )\Gamma
(2H+1 )^{1/({2H+1})}} {
(2H+1 )^2\sin({\uppi }/({2H+1}) )},
\]
where $\Gamma(\cdot)$ denotes the Gamma function.
\end{corollary}

Proofs of the statements are deferred to the \hyperref[app]{Appendix}.
Let us conclude the section by some comments on Theorems \ref{mainthm}
and \ref{thmexplicit}.

First, observe that (\ref{fisherintegral}) is of the type (\ref
{eqFInfogen}) for $\theta=\sigma^2$. The surprising fact is that one
can compute the integral $\int_0^\uppi  f^2(\lambda)/h_n^2(\lambda)
\mrmd \lambda$ obtaining expressions which for $\Diamond<4$ do not depend
on $f$ anymore. Let us shortly explain this. Suppose that Assumption
\ref{assX} holds with $\ell=1$. By classical results, we can conclude
that for $\lambda\downarrow0$, $f(\lambda) \sim C_\alpha\lambda
^{2\alpha}$ and $h_n(\lambda)\sim C_\alpha\sigma^2n^{-2\beta
}\lambda^{2\alpha}+\tau^2\lambda^{4K}$ with $C_\alpha=2\signtxt
(-\alpha)\Gamma(-2\alpha)\cos(\uppi \alpha)$. Next, observe that for
small $\lambda$, $f^2(\lambda)/h_n^2(\lambda)\approx n^{4\beta}$
whereas for large values the integrand behaves like $\lambda^{4\alpha
-4K}$. Now, $\Diamond\leq4$ is equivalent to $4\alpha-4K<-1$ and in
this case the integral will be determined in first order by
$f^2(\lambda)/h_n^2(\lambda)$ for small $\lambda$. Therefore, one
expects that $f$ and $h_n$ can be replaced by their corresponding small
value approximations $C_\alpha\lambda^{2\alpha}$ and $C_\alpha
\sigma^2n^{-2\beta}\lambda^{2\alpha}+\tau^2\lambda^{4K}$,
respectively, and
\[
\int_0^\uppi \frac{f^2(\lambda)}{h_n^2(\lambda)} \mrmd \lambda\approx
\sigma^{-4}n^{4\beta} \int_0^{\uppi }
\bigl(1+C_\alpha^{-1}\sigma^{-2}
\tau^2n^{2\beta}\lambda^{2K-2\alpha} \bigr)^{-2} \mrmd
\lambda.
\]
The r.h.s. can be explicitly solved and does not depend on $f$. In
contrast to that, for $\Diamond>4$, the Fisher information depends
also asymptotically on the whole spectrum $(0,\uppi ]$. This is why we
need the additional assumption (\ref{eqaddassumpifDiamondlarge})
which controls the continuity of $f$ globally.

The phase transition for $\Diamond=4$ does not only affect the
asymptotic constant but also leads to an elbow phenomenon in the rate
of convergence for estimation of $\sigma^2$. If $\Diamond\leq4$ the
optimal rate (neglecting sub-polynomial factors in the following) is
$n^{\Diamond\beta/2-1/2}$ whereas for $\Diamond>4$, it
is $n^{2\beta-1/2}$. The latter only depends on $\beta$.
Typically, if in an estimation problem an elbow effect occurs there are
different sources of errors which cannot be balanced and therefore the
best attainable rate is given by the maximum of the single error rates.
However, in our situation the optimal rate turns out to be the
\textit{minimum}, more precisely it is $\min(n^{\Diamond\beta/2
-1/2},n^{2\beta-1/2})$.

Let us also shortly remark on the dependence of the asymptotic Fisher
information on the scaling coefficient $\sigma$. Write $p=\Diamond/4$
and $\theta=\sigma^2$. Choosing $Q_n(\tau,\alpha,\Diamond)$
appropriately, we find that the Fisher information for $\theta$
observing $U\sim\mathcal{N}(\theta^p,Q_n)$ if $p<1$ and $U\sim
\mathcal{N}(\theta,Q_n)$ if $p\geq1$ coincides in first order with
(\ref{eqFisherInfoK-alphalarge})--(\ref
{eqFisherInfoK-alphamedium}) (standardizing $X$ such that $\sum
_{k=-\infty}^\infty\gamma_k^2=1$ if $p>1$). Hence, in an asymptotic
sense our original statistical estimation problem is related to a
Gaussian shift model where we want to estimate the $p$th power ($p\leq
1$) of the mean value.

Our results also cover the case $\alpha\in(-1/2,-1/4)$ for which the
autocovariance function is not (square) summable. In fact, the proof
turns out to be very subtle and requires quite restrictive conditions.
In particular, we have to impose an assumption on the increments of the
autocovariance which is much stronger than $\GM$.

In the critical case $\Diamond=4$ an additional log-factor appears in
the rate of convergence. In Theorem \ref{thmexplicit}, we have
restricted ourselves to the (most important) case where $\ell$ is a
power of the logarithm, which allows to evaluate the asymptotic Fisher
information in closed form. However, from the proof one can follow a
slightly more general version, namely that under the assumptions of
Theorem \ref{mainthm} (in particular $\log(n)\ell^2(n)\rightarrow
\infty$) and with $q_n:=n^{-4\beta}\ell^2(n^{4\beta})$,
\[
I_{\sigma^2}^n =n^{1-4\beta}\tau^{-4}\int
_{q_n}^1 \ell^2 \biggl(\frac1{
\lambda} \biggr)\frac{\dd\lambda}{\lambda} \bigl(1+\RMo(1)\bigr)+\RMo \bigl(nq_n
\log(n) \bigr).
\]

Theorems \ref{mainthm} and \ref{thmexplicit} are derived for
all $\alpha\in(-1/2,1/2), \alpha\neq0$. The case $\alpha=0$ is
indeed special since as $\alpha\rightarrow0$, $|C(\Diamond,\alpha
)|\rightarrow\infty$. However, in the fBM${}+{}$WN model (\ref
{eqmodfBM}) (recall that $H=1/2-\alpha$) this phenomenon does not
play a role because of $\ell(\lambda)=\mathrm{const.}=H |2H-1|$, which
converges to zero (fast enough) as $H\rightarrow1/2$. This explains
why $c_H$ is continuous for all $H\in(0,1)$, whereas in general
$\alpha\mapsto C(\Diamond, \alpha)$ is not.

Besides the classical case $H=1/2$, which was mentioned already in the
\hyperref[intro]{Introduction}, one can easily simplify the asymptotic Fisher information
in the fBM${}+{}$WN model (\ref{eqmodfBM}) for $H\in\{1/4, 3/4\}$.
Indeed, as a consequence of Corollary \ref{cortheoremholdsforfBM},
we obtain for the multiplicative inverse (which is the asymptotic
variance of our estimator)
\begin{eqnarray*}
H&=&1/4\mbox{:}\quad  \bigl(I_{\sigma^2}^n\bigr)^{-1}\sim
\frac{27}{\sqrt{3} \uppi
^{1/3}}\sigma^{8/3}\tau^{4/3}n^{-2/3}
\approx10.64 \sigma^{8/3}\tau^{4/3} n^{-2/3},
\\
H&=&1/2\mbox{:}\quad  \bigl(I_{\sigma^2}^n\bigr)^{-1}\sim8
\sigma^3\tau n^{-1/2},
\\
H&=&3/4\mbox{:}\quad  \bigl(I_{\sigma^2}^n\bigr)^{-1}\sim
\frac{25\sqrt{5+\sqrt
{5}}}{\sqrt{2} 3^{7/5} \uppi ^{1/5}}\sigma^{16/5}\tau^{4/5}n^{-2/5}
\approx8.12\sigma^{16/5}\tau^{4/5}n^{-2/5}.
\end{eqnarray*}
Finally, one should notice that the elbow effect observed in Gloter and
Jacod \cite{GloterJacod2001a} does not relate to our results. In fact
they have studied the fBM${}+{}$WN model for $H=1/2$ (i.e., BM${}+{}$WN), where the
variance of the noise is allowed to depend on $n$. With the notation of
model (\ref{eqmod}),~the change in the rate appears as $\beta
\downarrow0$. In particular, they also discuss the case $\beta< 0$ in
which the classical $n^{-1/2}$-rate can be achieved. In our framework,
$\beta<0$ corresponds to estimation of the~scaling parameter of $(Y_i)_i$.

\section{Examples}
\label{secexamples}

In the \hyperref[intro]{Introduction}, we have already discussed the main example of
estimating the scale parameter of fractional Brownian motion under
Gaussian measurement noise. The solution is given in (\ref
{eqFishInfo}) (cf. also Corollary \ref{cortheoremholdsforfBM}). In
order to provide some further illustration of the derived results, we
discuss two estimation problems for which the Fisher information can be
explicitly computed.

\textit{Large measurement error}: Let $(X_i)_i$ denote a stationary
process with long-range dependence. More precisely, assume that for
constants $A,C$, and self-similarity parameter $H\in(1/2,1)$, the
autocovariance satisfies $\gamma_k \sim Ak^{2H-2}$, $|\gamma
_{k+1}-\gamma_k|\leq C\gamma_k/k$ and for $H\in(1/2,3/4)$,
additionally $\sum_{k=-\infty}^\infty\gamma_k^2=1$. Suppose that we
observe the scaled process $(X_i)_i$ under large noise, that is,
\[
Z_{i,n}= \sigma X_i + \tau n^\beta
\varepsilon_{i},\qquad i=1,\ldots,n, 0<\beta<H-1/2.
\]
Now, with the notation of Theorem \ref{thmexplicit}, $\alpha=1/2-H$
and $\Diamond=2/(2H-1)$. In particular $\Diamond>4$ for $H\in
(1/2,3/4)$ and $\Diamond<4$ for $H\in(3/4,1)$. Therefore, an elbow
effect occurs at $H=3/4$ and the Fisher information is determined in
first order by
\begin{eqnarray*}
I_{\sigma^2}^n &\sim&\frac{n^{1-4\beta}}{2\tau^4},\qquad H\in\biggl(\frac12,
\frac34 \biggr),
\\
I_{\sigma^2}^n &\sim& 4\beta\frac{n^{1-4\beta}\log(n)}{\tau^4},\qquad H=\frac34,
\\
I_{\sigma^2}^n &\sim& n^{({2H-1-2\beta})/({2H-1})}\sigma^{-
({8H-6})/({2H-1})}\\
&&{}\times
\tau^{-2/({2H-1})} \frac{ (2A\Gamma(2H-1)\sin
(\uppi  H) )^{1/({2H-1})}}{(2H-1)^2\sin({\uppi }/({2H-1}))},\qquad H\in\biggl
(\frac34,1 \biggr).
\end{eqnarray*}

\textit{Integrated fractional Brownian motion}: Suppose that we are
interested in efficient estimation of $\sigma^2$ given observations
$(V_{1,n},\ldots,V_{n,n})$,
\[
V_{i,n}= \sigma\int_{i/n}^1
B_s^H\mrmd s +\tau\varepsilon_i,\qquad i=1,\ldots,n,
\varepsilon_i \stackrel{\mathrm{i.i.d.}} {\sim} \mathcal{N}(0,1),
H\in(0,1/4),
\]
and $(B_t^H)_t$ is a fBM, which is independent of the WN. With $\Delta
$ as in (\ref{backwarddiff}), $X_i:=\int_{i-1}^iB_s^H\mrmd s-\int_{0\vee
(i-2)}^{i-1} B_s^H\mrmd s$, and
$\bolds{\varepsilon}=(\varepsilon_1,\ldots,\varepsilon_n)$, consider
\[
\Delta\Delta^t (V_{1,n},\ldots,V_{n,n})^t
\stackrel{\mathcal{D}} {=}\bigl(\sigma n^{-1-H} X_i
\bigr)_{i=1,\ldots,n}+\tau\Delta\Delta^t\bolds{
\varepsilon}^t
\]
and note that $(X_i)_{i\geq2}$ is stationary. By defining $R_1$
appropriately, it is straightforward to verify Assumption \ref{assX}.
In particular, we find that (\ref{eqRfrobbdcond}) is bounded by
$\lesssim n^{4H-1}\leq1$ and $\gamma_k \sim H(2H-1)k^{2H-2}$.
Therefore, $\ell=H(1-2H)$, $\alpha=1/2-H, K=2$ and $\beta=1+H$ imply
\[
\Diamond=\frac{2}{2H+3}<4
\]
and the Fisher information is given by
\[
I_{\sigma^2}^n \sim n^{1/({2H+3})}\sigma^{-
({8H+10})/({2H+3})}
\tau^{-{2}/({2H+3})} \frac{(H+1)\sin^{1/({2H+3})} (\uppi  H )\Gamma
(2H+1 )^{1/({2H+3})}} {
(2H+3 )^2\sin({\uppi }/({2H+3}) )}.
\]

\begin{appendix}\label{app}
\section*{Appendix: Proofs}\label{secproofs}

\subsection{\texorpdfstring{Proof of Theorem \protect\ref{thmestisopt}}
{Proof of Theorem 1}}

\begin{lemmma}
\label{eqldsigmapreliminary} Let $\widetilde\sigma^2$ be as defined
in (\ref{eqprelimestdef}). Given $0<\sigmamin<\sigmamax<\infty$,
there exists an $N=N(\sigmamin, \sigmamax)$ such that for all $n\geq N$
and for all $\eta>0$,
\[
\sup_{\sigma\in[\sigmamin,\sigmamax]} \P\bigl(\bigl|\widetilde\sigma^2-
\sigma^2 \bigr|\geq\bigl(I_1^{A_n}
\bigr)^{-1/2} \max\bigl(1,\sigmamax^2\bigr)\eta\bigr)
\leq2\RMe ^{1/4-\eta/\sqrt{8}}.
\]
\end{lemmma}

\begin{pf}
Since $\delta_n \rightarrow0$, we can choose $N$ such that for all
$n\geq N$, $\delta_n \leq\sigmamin^2$ and $\delta_n^{-1}\geq
\sigmamax^2$. Then, for $\sigma^2 \in[\sigmamin, \sigmamax]$,
$|\widetilde\sigma^2-\sigma^2 | \leq|V-\sigma^2 |$. Let $\alpha
_{1},\alpha_2,\ldots$ be a sequence of i.i.d. $\chi_1^2$ random
variables. By Proposition 6 in Rohde and
D\"umbgen \cite{RohdeDuembgen2012} (similar statements have been
derived also elsewhere, for another reference see, for instance,
Johnstone \cite{Johnstone1999}, p. 74), for a vector $(\mu_i)_{i\in
A_n}$ of real-valued numbers
\[
\P\biggl( \biggl|\sum_{i\in A_n} \mu_i(
\alpha_i-1) \biggr| \geq\sqrt{2}\|\mu\|_2 \eta\biggr)
\leq2\RMe^{1/4-\eta/\sqrt{8}}.
\]
Note that in distribution, $V -\sigma^2
=
\sum_{i\in A_n} \mu_i (\alpha_i-1)$ with
\[
\mu_i= \bigl(2I_1^{A_n} \bigr)^{-1}
\frac{\lambda_i
n^{-2\beta}(\sigma^2n^{-2\beta}\lambda_i+1)}{(n^{-2\beta}\lambda_i+1)^2}.
\]
Application of the exponential
inequality above together with $\|\mu\|_2 \leq\break (2I_1^{A_n})^{-1/2}\max
(\sigma^2,1)$
yields the result.
\end{pf}

\begin{pf*}{Proof of Theorem \ref{thmestisopt}}
Due to the independence of $(\widetilde Z_i)_{i\in A_n}$ and
$(\widetilde Z_i)_{i\in A_n^c}$, the estimator $\widehat\sigma^2$ is
unbiased. In addition, using Lemma \ref{lempartitioneigvals}, the
theorem is proved once we have established that
\begin{eqnarray*}
\mbox{(I):}\hspace*{0.5pt}\quad \sup_{\sigma\in[\sigmamin, \sigmamax]} \frac{|{\Var}
(\widehat\sigma^2)-\E(I_{\widehat\sigma
^2}^{A_n^c})^{-1}|}{(I_{\sigma^2}^n)^{-1}}&=&\RMo(1)\quad\mbox{and}
\\
\mbox{(II):}\quad \sup_{\sigma\in[\sigmamin, \sigmamax]} \frac
{|\E(I_{\widehat
\sigma^2}^{A_n^c})^{-1}-(I_{\sigma^2}^{A_n^c})^{-1}|}{(I_{\sigma
^2}^n)^{-1}}&=&\RMo(1).
\end{eqnarray*}
In the following, we make frequently use of the following observation.
For any set
$B\subseteq\{1,\ldots,n\}$,
\[
\min\biggl(1, \biggl(\frac vu \biggr)^2 \biggr)I^B_v
\leq I_{u}^B\leq\max\biggl(1, \biggl(\frac vu
\biggr)^2 \biggr)I^B_v.
\]

(I): Writing $\E[\cdot]=\E[\E[\cdot|(\widetilde Z_{i,n})_{i
\in A_n}]]$,
\begin{eqnarray*}
\bigl|\Var\bigl(\widehat\sigma^2\bigr)-\E\bigl(I_{\widetilde\sigma
^2}^{A_n^c}
\bigr)^{-1} \bigr| &=& 2 \biggl|\E\biggl[\bigl(2I_{\widetilde\sigma^2}^{A_n^c}
\bigr)^{-2}\sum_{i\in
A_n^c}\frac{\lambda_i^3n^{-6\beta}(\sigma^2-\widetilde\sigma^2)}{
(\widetilde\sigma^2\lambda_in^{-2\beta}+1)^4}
\bigl(\bigl(\sigma^2+\widetilde\sigma^2\bigr)
\lambda_in^{-2\beta}+2 \bigr) \biggr] \biggr|
\\
&\leq& 2\E\biggl[\bigl(2I_{\widetilde\sigma^2}^{A_n^c}\bigr)^{-1} \biggl|
\frac
{\sigma^2}{\widetilde\sigma^2}-1 \biggr| \biggl(3+ \biggl|\frac{\sigma
^2}{\widetilde\sigma^2}-1 \biggr| \biggr) \biggr],
\end{eqnarray*}
where we used the inequalities
\[
\sum_{i \in A_n^c} \frac{\lambda_i^4n^{-8\beta}}{(\widetilde\sigma
^2 \lambda_in^{-2\beta}+1)^4} \leq\frac2{\widetilde
\sigma^4} I_{\widetilde\sigma^2}^{A_n^c} \quad\mbox{and}\quad \sum
_{i \in A_n^c} \frac{\lambda_i^3n^{-6\beta}}{(\widetilde\sigma
^2 \lambda_in^{-2\beta}+1)^4} \leq\frac1{\widetilde
\sigma^2} I_{\widetilde\sigma^2}^{A_n^c}.
\]
Therefore
\begin{eqnarray*}
\sup_{\sigma\in[\sigmamin, \sigmamax]} \E\biggl[\bigl(I_{\widetilde
\sigma^2}^{A_n^c}
\bigr)^{-1} \biggl|\frac{\sigma^2}{\widetilde\sigma
^2}-1 \biggr| \biggr] &=& \sup_{\sigma\in[\sigmamin, \sigmamax]}
\int\P\bigl[ \bigl|\sigma^2-\widetilde\sigma^2\bigr|\geq
\widetilde\sigma^2 I_{\widetilde\sigma
^2}^{A_n^c} t \bigr] \mrmd t
\\
&\leq& \int\sup_{\sigma\in[\sigmamin, \sigmamax]} \P\bigl[ \bigl|\sigma^2-
\widetilde\sigma^2\bigr|\geq\delta_n I_{1}^{A_n^c}
t \bigr] \mrmd t.
\end{eqnarray*}
Application of Lemma \ref{eqldsigmapreliminary} together with $\int
_0^\infty\exp(-At)\mrmd t= A^{-1}$ yields
\[
\sup_{\sigma\in[\sigmamin, \sigmamax]} \E\biggl[\bigl(I_{\widetilde
\sigma^2}^{A_n^c}
\bigr)^{-1} \biggl|\frac{\sigma^2}{\widetilde\sigma
^2}-1 \biggr| \biggr] \lesssim\bigl(
\delta_n^2 I_1^{A_n}
\bigr)^{-1/2} \bigl(I_1^{A_n^c}\bigr)^{-1}.
\]
Similar,
\begin{eqnarray*}
\sup_{\sigma\in[\sigmamin, \sigmamax]} \E\biggl[\bigl(I_{\widetilde
\sigma^2}^{A_n^c}
\bigr)^{-1} \biggl|\frac{\sigma^2}{\widetilde\sigma
^2}-1 \biggr|^2 \biggr] &=& \sup
_{\sigma\in[\sigmamin, \sigmamax]} \int\P\bigl[\bigl|\sigma^2-\widetilde
\sigma^2\bigr|\geq\widetilde\sigma^2 \bigl(
I_{\widetilde
\sigma^2}^{A_n^c}\bigr)^{1/2} \sqrt{t} \bigr] \mrmd t
\\
&\leq& \int\sup_{\sigma\in[\sigmamin, \sigmamax]} \P\bigl[\bigl|\sigma^2-
\widetilde\sigma^2\bigr|\geq\delta_n \bigl(
I_{1}^{A_n^c}\bigr)^{1/2} \sqrt{t} \bigr] \mrmd t
\end{eqnarray*}
and because of $\int_0^\infty\exp(-A\sqrt{t})\mrmd t= 2A^{-2}$,
\[
\sup_{\sigma\in[\sigmamin, \sigmamax]} \E\biggl[\bigl(I_{\widetilde
\sigma^2}^{A_n^c}
\bigr)^{-1} \biggl|\frac{\sigma^2}{\widetilde\sigma
^2}-1 \biggr|^2 \biggr] \lesssim
\bigl(\delta_n^2 I_1^{A_n}
\bigr)^{-1} \bigl(I_1^{A_n^c}\bigr)^{-1}.
\]
By definition $\delta_n^2 I_1^{A_n}\rightarrow\infty$. Since for
sufficiently large $n$, by Lemma \ref{lempartitioneigvals},
%
\begin{equation}\label{eqIsigma2bd}
I_{\sigma^2}^n\leq\max\bigl(1, \sigma^{-4}
\bigr)I_1^n \leq2\max\bigl(1, \sigma^{-4}
\bigr)I_1^{A_n^c},
\end{equation}
it follows that
\[
\sup_{\sigma\in[\sigmamin, \sigmamax]} \frac{|{\Var}(\widehat
\sigma^2)-\E
(I_{\widetilde\sigma^2}^{A_n^c})^{-1}|}{(I_{\sigma^2}^n)^{-1}} \lesssim
I_{1}^{A_n^c}
\sup_{\sigma\in[\sigmamin, \sigmamax]} \bigl|{\Var}\bigl(\widehat\sigma^2\bigr
)-\E
\bigl(I_{\widetilde\sigma
^2}^{A_n^c}\bigr)^{-1} \bigr| =\RMo (1 ).
\]

(II): From Taylor expansion, we find that for positive $x,y$,
$|x^{-2}-y^{-2}|\leq2(\min(x,y))^{-3} |x-y|$. Therefore, we can bound
\[
I_{\sigma^2}^{A_n^c}-I_{\widetilde\sigma^2}^{A_n^c} = \frac12 \sum
_{i \in A_n^c} n^{-4\beta}\lambda_i^2
\biggl(\frac1{(\sigma^2n^{-2\beta}\lambda_i+1
)^2}-\frac1{(\widetilde\sigma^2n^{-2\beta}
\lambda_i+1)^2} \biggr)
\]
by
\begin{eqnarray*}
\bigl|I_{\sigma^2}^{A_n^c}-I_{\widetilde\sigma^2}^{A_n^c} \bigr| &\leq& \sum
_{i \in A_n^c} \frac{n^{-6\beta}\lambda_i^3}{(\min(\sigma
^2,\widetilde\sigma^2)n^{-2\beta}\lambda_i+1)^3} \bigl|\sigma^2-
\widetilde\sigma^2\bigr|
\\
&=& 2 \max\bigl(\sigma^{-2}, \widetilde\sigma^{-2}\bigr)
I_{\sigma^2\wedge
\widetilde\sigma^2}^{A_n^c} \bigl|\sigma^2-\widetilde
\sigma^2\bigr| \\
&\leq&\bigl(\sigmamin^{-6}+1\bigr)
\delta_n^{-2} I_{\widetilde\sigma
^2}^{A_n^c} \bigl|
\sigma^2-\widetilde\sigma^2\bigr|.
\end{eqnarray*}
Thus,
\[
\sup_{\sigma\in[\sigmamin, \sigmamax]} I_{\sigma^2}^{A_n^c} \bigl| \E
\bigl(I_{\widetilde\sigma^2}^{A_n^c}\bigr)^{-1} -\bigl(I_{\sigma^2}^{A_n^c}
\bigr)^{-1} \bigr| \lesssim\delta_n^{-2} \int\sup
_{\sigma\in[\sigmamin, \sigmamax
]} \P\bigl(\bigl|\sigma^2-\widetilde
\sigma^2\bigr|\geq t\bigr)\mrmd t \lesssim\bigl(\delta_n^4
I_1^{A_n} \bigr)^{-1/2} \rightarrow0.
\]
Using (\ref{eqIsigma2bd}), the convergence in (II) follows and this
completes the proof of Theorem~\ref{thmestisopt}.
\end{pf*}

\subsection{\texorpdfstring{Notation and remarks for the results in Section \protect\ref{secclosedform}}
{Notation and remarks for the results in Section 2.2}}\label{notation}

Let us first give some notation. Whenever it is clear from the context,
we omit the index $n$. In particular, we suppress the index $n$ of the
spectral density $h_n=h$ and the estimator \mbox{$\widehat\sigma
^2_n=\widehat\sigma^2$}. Inequalities for Hermitian matrices should be
understand in the sense of partial Loewner ordering. The matrix norms
$\|\cdot\|_2$ and $\|\cdot\|_\infty$ denote the Frobenius and
spectral norm, respectively. Furthermore, we write $\wedge$ and $\vee
$ for the minimum/maximum and
%
\begin{equation}\label{equindef}
u_i:=u_{i,n}:=\uppi \frac{2i-1}{2n+1}.
\end{equation}
The last definition occurs frequently in connection with
finite-dimensional approximations due to the transformation property of
the discrete cosine transform. Let $\lfloor\cdot\rfloor$ be the
Gauss bracket. As in (\ref{eqdefrn}), we denote by $(r_n)_n$ the rate
at which the Fisher information tends to infinity (this still needs to
be proved, of course). For technical reasons, however, it will be
chosen in the following as an integer sequence, that is,
\[
r_n:= \cases{ \bigl\lfloor n^{1-{\beta}/({K-\alpha})} \bigl(\ell
\bigl(n^{{\beta
}/({K-\alpha})}\bigr) \bigr)^{1/({2K-2\alpha})}+n^{1-4\beta} \bigr\rfloor,
&\quad if $\Diamond\neq4$,
\vspace*{2pt}\cr
\bigl\lfloor n^{1-4\beta}\log(n)
\ell^2(n)\bigr\rfloor, &\quad if $\Diamond= 4$.}
\]
This definition will be used at many places throughout the proofs. In
particular, one should keep in mind that $u_{r_n}=\RMO(r_n/n)$ and as a
direct consequence of (\ref{assell}).
%
\begin{lemmma}
\label{lemelliter}
If $\ell$ is as in Assumption \ref{assX} and $\Diamond<4$, then
\[
\frac{\ell(u_{r_n}^{-1})}{\ell(n^{{\beta}/({K-\alpha})})} \rightarrow1.
\]
\end{lemmma}
The projection of a function on $\{\cos(k\cdot)\dvtx 0\leq k\leq n\}$ is
denoted by $S_n$, that is, if $g=\sum^\infty_{k=-\infty}g_k\cos
(k\cdot)$ with
$g_k=g_{-k}$, then
\[
S_ng= \sum_{k=-n}^n
g_k\cos(k\cdot).
\]
We write $T_n (g )$ for the $n\times n$ Toeplitz matrix
corresponding to
$g$, that is,
%
\begin{equation}
\label{eqdeftn}
T_n (g )=(g_{|i-j|})_{i,j=1,\ldots,n}.
\end{equation}
In particular, $\Cov(\bX')=T_n (f )$, with $f$ as in (\ref
{eqsdnotation}). Let $\dct$ be the discrete cosine transform
\begin{eqnarray*}
\dct &=& \biggl(\frac{2}{\sqrt{2n+1}}\cos\biggl[ \biggl(i-\frac
{1}{2} \biggr)
\biggl(j-\frac{1}{2} \biggr)\frac{2\uppi }{2n+1} \biggr] \biggr
)_{i,j=1,\ldots,n}\\
&=& \biggl(\frac{2}{\sqrt{2n+1}}\cos\biggl[ \biggl(i-\frac{1}{2}
\biggr)u_j \biggr] \biggr)_{i,j=1,\ldots,n}
\end{eqnarray*}
(which is DCT-VIII in the notation of \cite{BritanakYipRao2007}).
Note that $\dct=\dct^t$ is orthonormal. Further introduce the matrix
$D_n (g )$ as
%
\begin{equation}
\label{eqdefdn}D_n (g ):=\dct\cdot\operatorname{diag}
\bigl(g(u_1 ), g(u_2),\ldots,g(u_n) \bigr)
\cdot\dct.
\end{equation}
Asymptotically, the eigenvalues of $D_n (g )$ and $T_n
(g )$ are
`close', provided the symbol $g$ is sufficiently smooth (cf. Lemma C.4).

If $\Cov(\bY)=\tau^2(\Delta^t\Delta)^K$, consider the observation
vector $\widetilde{\mathbf{Z}}=(\widetilde Z_{1,n},\ldots,\widetilde
Z_{n,n})$ which is obtained by the
invertible transformation $\widetilde Z_{i,n}=Z_{n-i,n}$. These
observations satisfy our assumptions with
$\Cov(\bY)=\tau^2(\Delta\Delta^t)^K$. Therefore, without loss of
generality, we can and will consider only the first
case of Assumption \ref{assY}, that is, $\Cov(\bY)=\tau^2(\Delta\Delta
^t)^K$. By Lemma C.2, $\dct$ diagonalizes $\Delta^t \Delta$ and
hence also $\Cov(\bY)=\tau^2(\Delta^t \Delta)^K$. The
eigenvalues of $\Cov(\bY)$ are explicitly given by
$4^K\tau^2\sin^{2K}(u_i\frac{\uppi }2)$, $i=1,\ldots,n$.

For the subsequent proofs, the following three elementary inequalities
turn out to be very useful. Firstly, from (\ref
{eqautocorrasympbehavior}) and Potter's bound (cf., e.g.,
Bingham \textit{et al.}  \cite{BinghamGoldieTeugels1987}) it follows that
for any $\varepsilon>0$ there exists a $k_0$ such that for any $k\geq k_0$
%
\begin{equation}\label{eqfdmtlrXineqs}
\tfrac14 k^{-2\alpha-1-\varepsilon}\leq\tfrac12 k^{-2\alpha-1}\ell(k) \leq
\gamma_k \leq2 k^{-2\alpha-1}\ell(k) \leq4 k^{-2\alpha
-1+\varepsilon}.
\end{equation}
Moreover, we can find a constant $C_1$ such that $\gamma_k\leq C_1
k^{-2\alpha-1+\varepsilon}$ for all $k=1,2,\ldots$\,. Secondly, if
$f(\lambda) \sim Cn^{-2\beta} \lambda^{2\alpha}\ell(1/\lambda)$,
then, for every $\varepsilon>0$ we can find a $\delta>0$, such that for
all $\lambda\in(0,\delta]$,
%
\begin{equation}\label{eqfdmtlfXineqs}
\frac14 C \lambda^{2\alpha+\varepsilon} \leq\frac12 C \lambda^{2\alpha}
\ell
\biggl(\frac1{\lambda} \biggr) \leq f(\lambda) \leq2 C \lambda^{2\alpha
} \ell
\biggl(\frac1{\lambda} \biggr) \leq4 C \lambda^{2\alpha-\varepsilon}.
\end{equation}
Additionally, under the assumptions of Theorem \ref{mainthm}, we know
that $f$ is bounded on $[\delta, \uppi ]$ for every $\delta>0$ (cf.
Lemma C.6) and therefore the upper bound $f(\lambda) \leq4 C\lambda
^{2\alpha-\varepsilon}$ can be extended to all $\lambda\in(0,\uppi ]$ by
enlarging the constant appropriately. Finally, for all $\lambda\in
(0,\uppi ]$, we have
%
\begin{equation}\label{eqfdmtlhineqs}
\sigma^2n^{-2\beta}f(\lambda)+4^{-K}
\tau^2\lambda^{2K} \leq h(\lambda) \leq\sigma^2n^{-2\beta}f(
\lambda)+\tau^2\lambda^{2K}.
\end{equation}

\subsection{\texorpdfstring{Proof of Theorem \protect\ref{mainthm}}
{Proof of Theorem 2}}

The proof of the main theorem builds in a very neat way upon an
elementary analytical observation (cf. Lemma \ref
{ConvergenceofSequences}) which leads in a second step to a trace
approximation for positive semidefinite matrices (cf. Lemma \ref
{lemtracetrick}). This approximation result does not require any
assumption on the behavior of the smallest or largest eigenvalue.
Together with a rather standard but slightly technical Riemann
approximation argument, we can then deduce a generalized version of
Theorem \ref{mainthm}.

\begin{lemmma}\label{ConvergenceofSequences}
Let $(x_n)_n, (y_n)_n, (q_n)_n, (\omega_n)_n$ be positive sequences
such that $|\sqrt{x_n}-\sqrt{y_n}|=\RMO (q_n )$ and $\omega
_n\rightarrow\infty$. Then, $x_n=y_n(1+\RMo(1))+\RMO(q_n^2\omega_n)$.
\end{lemmma}
\begin{pf}
Set\vspace*{1pt} $a_n=\max(x_n,y_n)$ and $b_n=\min(x_n,y_n)$. We obtain
$(a_n/b_n)^{1/2}-1=\RMO(q_nb_n^{-1/2})$ implying
$a_n/b_n=1+\RMO(q_nb_n^{-1/2}+q_n^2 b_n^{-1})$ and
$a_n=b_n+\RMO(q_nb_n^{1/2}+q_n^2)$. Since $2q_nb_n^{1/2}\leq b_n/\omega
_n+q_n^2\omega_n$ and $\omega_n\rightarrow\infty$, we find
$a_n=b_n(1+\RMo(1))+\RMO(q_n^2\omega_n)$.
\end{pf}

Let $A$ be an $n\times n$ matrix. For convenience, we introduce the notation
\[
\langle A\rangle:= \operatorname{tr}\bigl(A^2\bigr).
\]

\begin{lemmma}
\label{lemtracetrick}
Let $A_1, A_2, B$ be (sequences of) positive semidefinite, $n\times n$
matrices and suppose that $A_1$ is invertible. If $(\omega_n)_n$ is a
positive sequence tending to infinity, then, for $n \rightarrow\infty$,
\[
\langle A_1B\rangle=\langle A_2B\rangle\bigl(1+\RMo(1)
\bigr) +\RMO \bigl(\bigl\langle B(A_1-A_2)\bigr\rangle
\omega_n \bigr).
\]
Furthermore, if $B\leq A_1^{-1}$, then
\[
\langle A_1B\rangle=\langle A_2B\rangle\bigl(1+\RMo(1)
\bigr) +\RMO \bigl(\bigl\langle\id_n-A_2^{1/2}A_1^{-1}A_2^{1/2}
\bigr\rangle\omega_n \bigr).
\]
\end{lemmma}

\begin{pf}
By Cauchy--Schwarz,
\begin{eqnarray*}
\bigl|\langle A_1B\rangle-\langle A_2B\rangle\bigr| &=& \bigl|
\operatorname{tr}\bigl(B(A_1-A_2)BA_1\bigr)+
\operatorname{tr}\bigl(BA_2B(A_1-A_2)\bigr)
\bigr|
\\
&\leq&\bigl\|B^{1/2}(A_1-A_2)B^{1/2}
\bigr\|_2 \bigl[\operatorname{tr}\bigl((A_1B)^2
\bigr)^{1/2}+\operatorname{tr}\bigl((A_2B)^2
\bigr)^{1/2} \bigr].
\end{eqnarray*}
For the last inequality we have rewritten $\operatorname
{tr}(B(A_1-A_2)BA_1)$ and
$\operatorname{tr}(BA_2B(A_1-A_2))$ as $\operatorname
{tr}(B^{1/2}(A_1-A_2)B^{1/2} \cdot B^{1/2}A_1B^{1/2})$ and
$\operatorname{tr}(B^{1/2}(A_1-A_2)B^{1/2} \cdot B^{1/2}A_2B^{1/2})$.
Since $\langle B(A_1-A_2)\rangle=\|
B^{1/2}(A_1-A_2)B^{1/2}\|_2^2$, the result follows with Lemma \ref
{ConvergenceofSequences}.

To prove the second claim, write
\[
B^{1/2}(A_1-A_2)B^{1/2}
=B^{1/2}A_1^{1/2}\bigl(\id_n-A_1^{-1/2}A_2A_1^{-1/2}
\bigr)A_1^{1/2}B^{1/2}
\]
and note that due to $A_1^{1/2}BA_1^{1/2}\leq\id_n$,
\[
\bigl\langle B(A_1-A_2)\bigr\rangle\leq\bigl\langle
\id_n-A_1^{-1/2}A_2A_1^{-1/2}
\bigr\rangle=\bigl\langle\id_n-A_2^{1/2}A_1^{-1}A_2^{1/2}
\bigr\rangle.
\]
%
\upqed\end{pf}


In the case $\alpha>0$ the multiplicative inverse of the spectral
density $h$ has a singularity at zero. In
order to deal with this, we introduce the regularized spectral density
$\widetilde h$, which is defined as follows: let
$(\rho_n)$ be a sequence of positive integers satisfying $\rho_n \ll
r_n$. Then, we define
\[
\widetilde h(\lambda):=\cases{ h(\lambda) \vee h(u_{\rho_n}), &\quad $\lambda
\leq u_{\rho_n}$,
\cr
h(\lambda), &\quad else,}
\]
with $u_{\rho_n}$ as in (\ref{equindef}). Replacing $h$ by $f$,
define in the same way $\widetilde f$. We will prove a generalized
version of Theorem \ref{mainthm} for a generic sequence $(\rho_n)_n$.
In a second step, the different versions of
the main theorem are deduced and $\rho_n$ will be chosen according to
the specific setting. Heuristically, we may interpret this spectral
regularization as adding an asymptotically noninformative (i.e.,
sufficiently small) WN process to our observation vector. This induces
some stability, which becomes important in the bounds for the inverse
covariance matrices.

\begin{theorem}
\label{thmgeneralizedmaintheorem}
Work under Assumptions \ref{assX} and \ref{assY} in model
(\ref{eqmod}). Suppose that $\alpha\in(-1/2,1/2)$ and $K-\alpha>
\beta\vee1/4$. If:
\begin{itemize}[(iii)]
\item[(i)] $(\gamma_k)_{k\geq0}$ is in $\GM$, $f$ is bounded on any
interval $[\delta,\uppi ]$ with $\delta>0$ and there exists a positive,
quasi-monotone slowly varying function $\ell$, such that
\[
f(\lambda) \sim2\signtxt(-\alpha)\Gamma(-2\alpha)\cos(\uppi \alpha)\lambda
^{2\alpha}\ell(1/\lambda),
\]
\item[(ii)] $n^{-4\beta-4\alpha-2+2\varepsilon}\sum_{i=1}^n (u_{i,n}
\widetilde h(u_{i,n}))^{-2} =\RMo(r_n)$, for some $\varepsilon>0$,
\item[(iii)] $ \langle D_n^{-1} ({\widetilde h})(D_n
(S_n f )-T_n (f )) \rangle+ \sup_{\lambda\in(0,\uppi
]} \widetilde
h^{-2}(\lambda) =\RMo(r_n n^{4\beta})$.
\end{itemize}
Then, the asymptotic Fisher information of $\sigma^2$ is
%
\begin{equation}\label{fisherintegralappendix}
I_{\sigma^2}^n = \frac{n^{1-4\beta}}{2\uppi }\int^\uppi _0
\frac
{f^2(\lambda)}{h^2(\lambda)} \mrmd \lambda\bigl(1+\RMo(1)\bigr)+\RMo(r_n).
\end{equation}
If the condition $K-\alpha> \beta\vee1/4$ is replaced by the weaker
assumption $K-\alpha>\beta$, imposing additionally $\log(n)\ell
^2(n)\rightarrow\infty$ in the critical case $K-\alpha=1/4$, then
\textup{(\ref{fisherintegralappendix})} holds, provided there exists a
constant $c_f$ such that
%
\begin{equation}\label{eqfX-fYestlargelambda}
\bigl|f(\lambda)-f(\mu) \bigr|\leq c_f n^{-2\beta}
\lambda^{2\alpha
-2}|\lambda-\mu|\qquad \mbox{for all } 0<\lambda\leq\mu\leq\uppi .
\end{equation}
\end{theorem}

\begin{remark}
\label{remondiffofgenmainthm}
Later on we will see that the different parts of Theorem \ref{mainthm}
follow from Theorem \ref{thmgeneralizedmaintheorem}. If $(X_i')_i$
has long-memory, condition (iii) turns out to be quite difficult to
verify. Although it would be easier (and more standard) to formulate
the condition with respect to squared Frobenius norms, let us shortly
explain, why the use of the $\langle\cdot\rangle$ notation is
essential. By definition, $\langle A\rangle=\operatorname{tr}(A^2)$
for an
$n\times n$, square matrix $A$, which in turn can be upper bounded by
the squared Frobenius norm of $A$ (cf. Lemma C.1(i)). This is even an
identity if $A$ is symmetric but can be very rough in general. To see
this consider $(A)_{i,j}=1/i$. Then $\langle A\rangle=\log^2n$ but
the squared Frobenius norm is of order $n$ which is much worse. Since
these phenomenons occur in some cases, condition (iii) is stated for
squared traces.
\end{remark}

\begin{pf*}{Proof of Theorem \ref{thmgeneralizedmaintheorem}}
Recall the explicit expression for the Fisher information in (\ref
{eqFishexplexpr}). The proof is subdivided into three steps, namely
\begin{eqnarray*}
\bigl\langle D_n (\widetilde f )D_n^{-1} ({
\widetilde h}) \bigr\rangle&=&\frac{n}{\uppi }\int^\uppi _0
\frac{f^2(\lambda)}{h^2(\lambda)}\mrmd  \lambda\bigl(1+\RMo(1) \bigr)+\RMo \bigl(r_n
n^{4\beta} \bigr),
\\
\bigl\langle T_n (f )D_n^{-1} ({\widetilde h})
\bigr\rangle&=& \bigl\langle D_n (\widetilde f )D_n^{-1}
({\widetilde h}) \bigr\rangle\bigl(1+\RMo(1) \bigr)+\RMo
\bigl(r_n n^{4\beta} \bigr),
\\
2 I_{\sigma^2}^n= \bigl\langle n^{-2\beta}T_n
(f )\Cov(\bZ)^{-1} \bigr\rangle&=& \bigl\langle n^{-2\beta}T_n
(f )D_n^{-1} ({\widetilde h}) \bigr\rangle\bigl(1+\RMo(1)
\bigr)+\RMo (r_n ),
\end{eqnarray*}
which are denoted by (I), (II) and (III), respectively.

(I): By the trivial bound $\int_0^\uppi  f^2(\lambda)/h^2(\lambda
) \mrmd \lambda=\RMO(n^{4\beta})=\RMo(r_nn^{4\beta})$, we can replace the
normalization factor $n/\uppi $ by $(2n+1)/(2\uppi )$. Thus, using (\ref
{eqdefdn}), it is sufficient to show that
%
\begin{equation}\label{eqredstep1sumdiff}
\Biggl|\sum_{i=1}^n \frac{\widetilde f^2(u_i)}{\widetilde h^2(u_i)} -
\frac{2n+1}{2\uppi } \int_0^\uppi
\frac{f^2(\lambda)}{h^2(\lambda)} \mrmd \lambda\Biggr|=\RMo\bigl(r_n n^{4\beta}\bigr).
\end{equation}
Now, let us treat the cases $K-\alpha>1/4 \vee\beta$ and $\beta
<K-\alpha\leq1/4$, separately.

\textit{If $K-\alpha>1/4 \vee\beta$ holds}: Using assumption (i) and
$r_n\ll n$, we can find integer sequences $(r_n^+)$ and $(r_n^-)$ such
that $\rho_n \ll r_n^-\ll r_n\ll r_n^+\ll n$ and
$(q_n+(r_n^-)^{-1})r_n^+=\RMo(r_n)$ with
%
\begin{equation}\label{eqqndef}
q_n:=q_n\bigl(r_n^+\bigr):=\mathop{\sup
_{0<\lambda\leq u_{r_n^+}}}_{1\leq\mu
\leq2} \biggl| \frac{C\lambda
^{2\alpha}\ell(1/\lambda)}{f(\lambda)}-1 \biggr| + \biggl| 1-
\frac{\ell({\mu}/{\lambda} )}{\ell
(1/{\lambda} )} \biggr|.
\end{equation}
Since $\sigma^2n^{-2\beta} f \leq h$ and $\sigma^2 n^{-2\beta}
\widetilde f\leq
\widetilde h$, it follows that
\[
\sum_{i=1}^{r_n^-} \frac{\widetilde f^2(u_i)}{\widetilde h^2(u_i)} =\RMo
\bigl(r_n n^{4\beta}\bigr) \quad\mbox{and}\quad n \int
_0^{u_{r_n^{-}}} \frac{f^2(\lambda)}{h^2(\lambda)}\mrmd  \lambda=\RMo
\bigl(r_n n^{4\beta}\bigr)
\]
and together with Proposition C.1, we see that in (\ref
{eqredstep1sumdiff}) the sum over $i=1,\ldots,r_n^-$ and
$i=r_n^+,\ldots,n$ as well as the integral over $(0,u_{r_n^-}]\cup
[u_{r_n^+},\uppi ]$ are of order $\RMo(r_n)$ and thus negligible. Thus, we
have proved (I), once we have verified that
%
\begin{equation}\label{eqtoshowforI}
\Biggl|\sum_{i=r_n^-+1}^{r_n^+} \frac{f^2(u_i)}{h^2(u_i)} -
\frac{2n+1}{2\uppi } \int_{u_{r_n^-}}^{u_{r_n^+}}
\frac{f^2(\lambda
)}{h^2(\lambda)} \mrmd \lambda\Biggr|=\RMo\bigl(r_n n^{4\beta}\bigr).
\end{equation}
To see this, write
%
\begin{eqnarray}\label{eqreductionofintegral}
&&
\Biggl|\sum_{i=r_n^-+1}^{r_n^+} \frac{f^2(u_i)}{h^2(u_i)} -
\frac{2n+1}{2\uppi } \int_{u_{r_n^-}}^{u_{r_n^+}}
\frac{f^2(\lambda
)}{h^2(\lambda)} \mrmd \lambda\Biggr|\nonumber\\
&&\quad\leq \sum_{i=r_n^-+1}^{r_n^+}
\sup_{\xi_i\in[u_{i-1},u_i]} \biggl|\frac
{f^2(u_i)}{h^2(u_i)} -\frac{f^2(\xi_i)}{h^2(\xi_i)} \biggr|
\\
&&\quad\leq \frac{2n^{4\beta}}{\sigma^2} \sum_{i=r_n^-+1}^{r_n^+}
\sup_{\xi_i\in[u_{i-1},u_i]} \biggl|\frac
{n^{-2\beta}f(u_i)}{h(u_i)} -\frac{n^{-2\beta}f(\xi_i)}{h(\xi_i)}
\biggr|.\nonumber
\end{eqnarray}
Fix\vspace*{1pt} $i\in\{r_n^{-}+1,\ldots,r_n^+\}$ and let $a_1=\sigma^2n^{-2\beta
} f(u_i), a_2=\sigma^2n^{-2\beta} f(\xi_i), b_1=4^K\tau^2\sin
^{2K}(u_i/2)$, $b_2=4^K\tau^2\sin^{2K}(\xi_i/2), c_1=h(u_i),
c_2=h(\xi_i)$. Since $0\leq a_1\leq c_1=a_1+b_1$ and $0\leq a_2\leq
c_2=a_2+b_2$, we find
%
\begin{equation}\label{eqelemdiffform}
\biggl|\frac{a_1}{c_1}-\frac{a_2}{c_2} \biggr| \leq\frac{|a_1-a_2|+|b_1-b_2|}{c_1
\vee c_2}.
\end{equation}
Thus, for sufficiently large $n$,
%
\begin{eqnarray}\label{eqa1-a2diffest}
|a_1-a_2| &\leq& (a_1+a_2)q_n
+ \sigma^2Cn^{-2\beta} \biggl| u_i^{2\alpha}\ell
\biggl(\frac1{u_i} \biggr) -\xi_i^{2\alpha}\ell
\biggl(\frac1{\xi_i} \biggr) \biggr|
\nonumber
\\
&\leq& (3a_1+a_2)q_n + \sigma^2Cn^{-2\beta}
\ell\biggl(\frac1{\xi_i} \biggr) \bigl|u_i^{2\alpha}-
\xi_i^{2\alpha}\bigr|
\nonumber\\[-8pt]\\[-8pt]
&\leq& (3a_1+a_2)q_n + \frac{\uppi }n
\sigma^2Cn^{-2\beta} \ell\biggl(\frac1{\xi_i}
\biggr)\xi_i^{2\alpha-1} \nonumber\\
&\leq&(3a_1+a_2)q_n
+ 6\bigl(r_n^-\bigr)^{-1} a_2.\nonumber
\end{eqnarray}
On the other hand, we find $|b_1-b_2|\leq\frac{2K\uppi }{2n+1}
4^K\tau^2 \sin^{2K-1} (\frac{u_i}{2} )\leq
8K(r_n^{-})^{-1}b_1$, and therefore
\[
\biggl|\frac{a_1}{c_1}-\frac{a_2}{c_2} \biggr| \leq4q_n+(6+8K)
\bigl(r_n^-\bigr)^{-1}.
\]
Due to $(q_n+(r_n^-)^{-1})r_n^+=\RMo(r_n)$, we see by (\ref
{eqreductionofintegral}) that (\ref{eqtoshowforI}) is bounded by
$\RMo(r_n n^{4\beta})$. This completes the proof for part (I) if
$K-\alpha>1/4 \vee\beta$.

\textit{If $\beta<K-\alpha\leq1/4$}: The proof is very similar to the
one for the first case. Note that the assumptions imply $r_nn^{4\beta
}\gtrsim n$, $K=0$, and $\alpha\in[-1/4,-\beta)$. Similar as above
we see that it is sufficient to prove (\ref{eqtoshowforI}) for
$r_n^+=n$ and any sequence $r_n^-=\RMo(r_n)$. We may assume that
$r_n^-\rightarrow\infty$. Since by (\ref
{eqfX-fYestlargelambda}), $f$ and $h$ are continuous, we can apply
the mean value theorem, that is, for any $i$ there is a $\xi_i \in
(u_{i-1},u_i]$ with
%
\begin{equation}\label{eqmeanvaluefXh}
\frac{f^2(\xi_i)}{h^2(\xi_i)} = \frac{2n+1}{2\uppi } \int_{u_{i-1}}^{u_i}
\frac{f^2(\lambda)}{h^2(\lambda)} \mrmd \lambda
\end{equation}
and
%
\begin{eqnarray}\label{eqmeanvalueapplied}
&&\Biggl|\sum_{i=r_n^-+1}^{n} \frac{f^2(u_i)}{h^2(u_i)} -
\frac{2n+1}{2\uppi } \int_{u_{r_n^-}}^{u_n}
\frac{f^2(\lambda
)}{h^2(\lambda)} \mrmd \lambda\Biggr|\nonumber\\[-8pt]\\[-8pt]
&&\quad \leq\sum_{i=r_n^-+1}^{n}
\biggl(\frac{f(u_i)}{h(u_i)}+\frac{f(\xi
_i)}{h(\xi_i)} \biggr) \biggl|\frac{f(u_i)-f(\xi_i)}{h(u_i) \vee h(\xi
_i)} \biggr|.\nonumber
\end{eqnarray}
Pick an integer sequence $w_n$ such that $w_n=\RMo(n)$ and for some
$\varepsilon>0$, $n^{-4\alpha+\varepsilon}w_n^{4\alpha-\varepsilon
-1}=\RMo(1)$. Let $q_n(w_n)$ be as in (\ref{eqqndef}) with $r_n^+$
replaced by $w_n$. Note that since $u_{w_n}\rightarrow0$, the sequence
$(q_n(w_n))_n$ tends to zero. As in (\ref{eqa1-a2diffest}), we find
for sufficiently large $n$ and for all $i=r_n^-,\ldots,w_n$,
%
\begin{equation}\label{eqfXdiffsmallvals}
\bigl|f(u_i)-f(\xi_i)\bigr|\leq\bigl(3 f(u_i)+f(
\xi_i)\bigr)q_n(w_n)+2\sigma^2\uppi
Cn^{-1}\xi_i^{2\alpha-1-\varepsilon}.
\end{equation}
For large indices, we use the estimate (\ref
{eqfX-fYestlargelambda}), that is, $|f(u_i)-f(\xi_i)|\leq c_f \uppi
n^{-1}\xi_i^{2\alpha-2}$ for $i=w_n+1,\ldots,n$. Split the sum in
(\ref{eqmeanvalueapplied}) into $\sum_{i=r_n^-+1}^{w_n}+\sum
_{i=w_n+1}^n$. It is easy to bound the second sum, using the definition
of $(w_n)_n$, which in turn implies that $n^{-1} \sum_{i=w_n+1}^n \xi
_i^{4\alpha-2-\varepsilon}\lesssim n^{1-4\alpha+\varepsilon}w_n^{4\alpha
-1-\varepsilon}=\RMo(n)$. Similar, computing the sum $\sum
_{i=r_n^-+1}^{w_n}$, over the second term in (\ref
{eqfXdiffsmallvals}) yields for small $\varepsilon$, $n^{-1} \sum
_{i=r_n^-+1}^{w_n} \xi_i^{4\alpha-1-2\varepsilon}
\lesssim n^{-4\alpha+2\varepsilon}w_n^{4\alpha-2\varepsilon}=\RMo(n)$. For
the first term in (\ref{eqfXdiffsmallvals}), we see that (\ref
{eqmeanvalueapplied}) can be\vspace*{1pt} further bounded by a multiple of
$q_n(w_n)(\sum_{i=1}^n f^2(u_i)/h^2(u_i)+(2n+1)(2\uppi )^{-1}\int_0^\uppi
f^2(\lambda)/h^2(\lambda) \mrmd \lambda)$. Putting all estimates
together, we have derived
\[
\Biggl|\sum_{i=r_n^-+1}^{n} \frac{f^2(u_i)}{h^2(u_i)} -
\frac{2n+1}{2\uppi } \int_{u_{r_n^-}}^{u_n}
\frac{f^2(\lambda
)}{h^2(\lambda)} \mrmd \lambda\Biggr| \lesssim\frac{2n+1}{2\uppi } \int
_{0}^{\uppi } \frac{f^2(\lambda
)}{h^2(\lambda)} \mrmd \lambda+\RMo(n).
\]
This finishes the proof of part (I).

To prove (II) and (III) it will not be necessary to distinguish
whether $K-\alpha>1/4$ or $K-\alpha\leq1/4$.

(II): In Lemma \ref{lemtracetrick}, set $A_1=T_n(f)$,
$A_2=D_n(\widetilde f)$ and $B=D_n^{-1}(\widetilde h)$. We have to show that
%
\begin{equation}\label{eqtoshowforIIingenmt}
\bigl\langle B(A_1-A_2)\bigr\rangle= \bigl\langle
D_n^{-1}(\widetilde h) \bigl(T_n(f)-D_n(
\widetilde f) \bigr) \bigr\rangle= \RMo\bigl(r_n n^{4\beta}\bigr).
\end{equation}
First, note that due to (\ref{eqdefdn}), $\sigma^2n^{-2\beta}f\leq
\sigma^2n^{-2\beta}\widetilde f\leq\widetilde h$, and $\rho_n\ll r_n$,
\[
\bigl\langle D_n^{-1} ({\widetilde h})
\bigl(D_n (f )-D_n (\widetilde f ) \bigr) \bigr\rangle=
\sum_{j=1}^{\rho_n} \frac{(f(u_j)-\widetilde f(u_j))^2}{(\widetilde
h(u_j))^2}\leq2
\frac{\rho_n n^{4\beta}}{\sigma^4}=\RMo\bigl(r_n n^{4\beta}\bigr).
\]
By assumption $(\gamma_k)_k\in\GM$ and $\alpha\in(-1/2,1/2)$.
Therefore, we can use the estimate from Lemma C.3(i) together with
(\ref{eqfdmtlrXineqs}), that is, there exists a constant $C_1$,
such that
\[
\bigl|f(x)-S_nf(x)\bigr|\leq C_1 \frac1x \Biggl(|
\gamma_n|+\sum_{k=n+1}^\infty
\frac{|\gamma_k|}{k} \Biggr) \lesssim\frac1x n^{-2\alpha
-1+2\varepsilon} \qquad\mbox{for all }
x\in\biggl[\frac1n, \uppi \biggr),
\]
where the second inequality holds for sufficiently large $n$ and
$\varepsilon$ small. With (\ref{eqdefdn}) and assumption (ii) this yields
%
\begin{equation}\label{eqDnSnfXapprox}
\bigl\langle D_n^{-1}(\widetilde h) \bigl(D_n(S_n
f)-D_n(f) \bigr) \bigr\rangle= \RMo\bigl(r_n
n^{4\beta}\bigr).
\end{equation}
Decompose $T_n(f)-D_n(\widetilde f)= (T_n(f)-D_n(S_nf))+
(D_n(S_nf)-D_n(f))+(D_n(f)-D_n(\widetilde f))$. Now by Lemma C.1, (iv)
and assumption (iii), (\ref{eqtoshowforIIingenmt}) follows.

(III): Set $A_1=\Cov(\bZ)^{-1}, A_2=D_n^{-1}(\widetilde h)$ and
$B=\sigma^2n^{-2\beta} \Cov(\bX)$ and apply Lemma \ref
{lemtracetrick}. Since $B\leq\Cov(\bZ)=A_1^{-1}$ it is sufficient
to show
%
\begin{equation}\label{eqtoproveinstep3}
\bigl\langle\id_n -A_2^{1/2}A_1^{-1}A_2^{1/2}
\bigr\rangle= \bigl\langle\id_n -D_n^{-1/2}(
\widetilde h)\Cov(\bZ) D_n^{-1/2}(\widetilde h) \bigr\rangle
=\RMo(r_n).
\end{equation}
By the perfect diagonalization property of $\Cov(\bY)$ (cf. the
remarks in Section \ref{notation}), we have $\Cov(\bY)=D_n(h-\sigma
^2n^{-2\beta}f)$ and
\begin{eqnarray*}
\Cov(\bZ)&=& \sigma^2n^{-2\beta} \bigl[\Cov(\bR)+\Cov\bigl(
\bX',\bR\bigr)+\Cov\bigl(\bR,\bX'\bigr) \bigr]\\
&&{} +
\sigma^2n^{-2\beta} \bigl(T_n(f)-D_n(f)
\bigr)+D_n(h).
\end{eqnarray*}
Together with Lemma C.1(iv),
%
\begin{eqnarray}\label{eqDNineqIII}
&&
\bigl\langle\id_n -D_n^{-1/2}(\widetilde h)\Cov(
\bZ) D_n^{-1/2}(\widetilde h) \bigr\rangle \nonumber\\
&&\quad\lesssim \bigl
\langle\id_n -D_n\bigl(\widetilde h^{-1} h\bigr)
\bigr\rangle
\nonumber
\\
&&\qquad{}+n^{-4\beta}\sup_{\lambda\in(0,\uppi ]} \widetilde h^{-2}(
\lambda) \bigl\|{\Cov}(\bR)+\Cov\bigl(\bX',\bR\bigr)+\Cov\bigl(\bR,
\bX'\bigr) \bigr\|_2^2
\\
&&\qquad{}+ n^{-4\beta} \bigl\langle D_n^{-1}(\widetilde h)
\bigl(D_n(f)-D_n(S_nf) \bigr) \bigr\rangle
\nonumber
\\
&&\qquad{}+ n^{-4\beta} \bigl\langle D_n^{-1}(\widetilde h)
\bigl(D_n(S_nf)-T_n(f) \bigr) \bigr\rangle.\nonumber
\end{eqnarray}
For the first term note that because of (\ref{eqdefdn}) and $0\leq
1-h/\widetilde h\leq1$,
\[
\bigl\langle\id_n -D_n\bigl(\widetilde h^{-1}
h\bigr) \bigr\rangle=\sum^{\rho
_n}_{i=1}
\biggl(1-\frac{h(u_i)}{\widetilde{h}(u_i)} \biggr)^2=\RMo(r_n).
\]
The other three terms on the r.h.s. of (\ref{eqDNineqIII}) can be
seen to be of order $\RMo(r_n)$ as well, by Assumption \ref{assX},
assumption (iii), (\ref{eqDnSnfXapprox}), and assumption (iii)
again. Therefore, (\ref{eqtoproveinstep3}) holds and the proof is
completed.
\end{pf*}

\begin{pf*}{Proof of Theorem \ref{mainthm}, part \ref{version1}}
We check the conditions of Theorem \ref{thmgeneralizedmaintheorem}.
Note that the special case, that is, $K-\alpha\leq1/4$ implies
together with $K-\alpha>\beta$ that $K=0$ and $\alpha\in
[-1/4,-\beta)$. Hence, this case only plays a role in parts \ref
{version2} and \ref{version3} (in the latter only if $K=0$ and
$\alpha=-1/4$). All the derived estimates will work for both
situations and thus, in the following, we do not distinguish between
these two cases explicitly.

Let
$\rho_n=n^{1-(4\alpha)^{-1}(1-{\beta}/({K-\alpha}))+\delta}
\ll r_n$ for some $\delta>0$. Such a $\delta$ always exists thanks to
the assumption $K-\alpha>(4\alpha+1)\beta$. This assures that
%
\begin{equation}
(\rho_n/n)^{-4\alpha-2\varepsilon}=\RMo(r_n) \qquad\mbox{for $\varepsilon$
small enough.}
\end{equation}

\mbox{}\hphantom{ii}(i): By \cite{Tikhonov2007}, $(\gamma_k)_k \in\GM$. The
second part follows from Lemma C.6.

\mbox{}\hphantom{i}(ii): Making use of inequalities (\ref{eqfdmtlfXineqs}) and
(\ref{eqfdmtlhineqs}),
\begin{eqnarray*}
&&
n^{-4\beta-4\alpha-2+2\varepsilon}\sum_{i=1}^n
\bigl(u_{i,n} \widetilde h(u_{i,n})\bigr)^{-2}
\\
&&\!\!\quad\lesssim n^{-4\beta-4\alpha-2+2\varepsilon} \Biggl[\sum^{\rho
_n}_{i=1}
\biggl(\frac ni \biggl(\frac{n}{\rho_n} \biggr)^{2\alpha+\varepsilon
}n^{2\beta
}
\biggr)^2+\sum^{r_n}_{i=\rho_n+1}
\biggl(\frac{n}{i} \biggr)^{4\alpha
+2\varepsilon+2}n^{4\beta
}+\sum
^{n}_{i=r_n+1} \biggl(\frac{n}{i}
\biggr)^{4K+2} \Biggr]
\\
&&\!\!\quad\lesssim n^{4\varepsilon}(\rho_n)^{-4\alpha-2\varepsilon}+n^{4(K-\alpha
-\beta)+2\varepsilon}r_n^{-4K-1}=\RMo(r_n),
\end{eqnarray*}
if $\varepsilon$ is chosen small enough.

(iii): By Lemma C.1(iii) and (i), Lemma C.4, and $T_n (f
)=T_n (S_nf )$,
\begin{eqnarray*}
\bigl\langle D_n^{-1} ({\widetilde h})
\bigl(D_n (S_nf )-T_n (f )\bigr) \bigr\rangle
&\leq&\bigl\|D_n^{-1} ({\widetilde h})\bigr\|^2_\infty
\bigl\langle D_n (S_nf )-T_n (f )\bigr\rangle
\\
&\leq&\bigl\|D_n^{-1} ({\widetilde h})\bigr\|^2_\infty
\bigl\|D_n (S_nf )-T_n (f )\bigr\|_2^2\\
&\lesssim&\biggl(\frac{\rho
_n}n \biggr)^{-4\alpha-2\varepsilon}n^{4\beta}=\RMo
\bigl(r_n n^{4\beta}\bigr),
\end{eqnarray*}
if $\delta$ and $\varepsilon$ are chosen appropriately. By the same
arguments $\sup_{\lambda\in(0,\uppi ]} \widetilde h^{-2}(\lambda
)=\RMo(r_nn^{4\beta})$ and this completes the proof of the claim.
\end{pf*}

\begin{pf*}{Proof of Theorem \ref{mainthm}, part \ref{version2}}
Let $\rho_n=1$, that is, $\widetilde h(u_i)=h(u_i)$. The proof of this
part is similar to the one for (i). Again, we check the conditions of
Theorem \ref{thmgeneralizedmaintheorem}:

\mbox{}\hphantom{ii}(i): This follows from Lemma C.6.

\mbox{}\hphantom{i}(ii): Splitting the sum $\sum_{i=1}^n = \sum_{i=1}^{r_n}+\sum
_{i=r_n+1}^n$, we find for small $\varepsilon$, by inequalities (\ref
{eqfdmtlfXineqs}) and (\ref{eqfdmtlhineqs}),
\[
n^{-4\beta-4\alpha-2+2\varepsilon}\sum_{i=1}^n
\bigl(u_{i,n} \widetilde h(u_{i,n})\bigr)^{-2}
\lesssim n^{4\varepsilon}+n^{4(K-\alpha-\beta)+2\varepsilon}r_n^{-4K-1}=\RMo(r_n).
\]

(iii): Observe that by (\ref{eqfdmtlfXineqs}) and (\ref
{eqfdmtlhineqs}), $h(\lambda)\gtrsim n^{-{2K\beta}/({K-\alpha
-\varepsilon/2})}$. Together with Lemma C.4,
\begin{eqnarray*}
\bigl\langle D_n^{-1} ({\widetilde h})
\bigl(D_n (S_nf )-T_n (f )\bigr) \bigr\rangle
+\sup_{\lambda\in(0,\uppi ]} \widetilde h^{-2}(\lambda) &\leq& \sup
_{\lambda\in(0,\uppi ]} \widetilde h^{-2}(\lambda) \bigl( \bigl\|
D_n (S_nf )-T_n (f ) \bigr\|_2^2+1
\bigr)
\\
&\lesssim& n^{{4K\beta}/({K-\alpha-\varepsilon/2})-4\alpha+2\varepsilon
}=\RMo\bigl(r_nn^{4\beta}\bigr).
\end{eqnarray*}
\upqed\end{pf*}

\begin{pf*}{Proof of Theorem \ref{mainthm}, part \ref{version3}}
We apply Theorem \ref{thmgeneralizedmaintheorem} with $\widetilde
f(u_i)=f(u_i)$ and $\widetilde h(u_i)=h(u_i)$, that is, $\rho_n=1$.

\mbox{}\hphantom{ii}(i): This follows from Lemmas C.5 and C.6.

For the following parts we make frequently use of the inequalities
(\ref{eqfdmtlrXineqs})--(\ref{eqfdmtlhineqs}) and the
subsequent comments.

\mbox{}\hphantom{i}(ii): Since $u_{r_n}\rightarrow0$, we find, if $n$ is
sufficiently large,
\[
\sum_{i=1}^n \frac1{u_{i}^2
h(u_i)^2} \lesssim\sum_{i=1}^{r_n}
\frac{n^{4\beta}}{u_{i}^2 f(u_i)^2} + \sum_{i=r_n+1}^n
\frac1{u_i^{2+4K}} \lesssim n^{2+4\beta+4\alpha+2\varepsilon
}r_n^{-1-4\alpha}+
n^{2+4K}r_n^{-1-4K}.
\]
Therefore,
\[
n^{-4\beta-4\alpha-2+2\varepsilon} \sum_{i=1}^n
\frac1{u_{i}^2 h(u_i)^2} \lesssim
n^{4\varepsilon}r_n^{-1-4\alpha}+n^{-4\beta-4\alpha
+4K+2\varepsilon}r_n^{-1-4K}
=\RMo(r_n)
\]
for $\varepsilon$ sufficiently small.

(iii): It is straightforward to bound $\sup_{\lambda\in(0,\uppi
]} \widetilde h^{-2}(\lambda)$ by a multiple of $n^{{4K\beta
}/({K-\alpha-\varepsilon/2})}=\RMo(r_nn^{4\beta})$, which immediately implies
that the second term has the right order. However, to show the same
rate for the first term turns out to be the most difficult part of the
proof. Let us shortly remark on that. The crucial point is that
although we have good control on the spectral density~$h$, this gives
no direct link to entries of the inverse of $D_n(h)$. In contrast to
the proofs above, estimating $\langle
D_n^{-1}(h)(D_n(S_nf)-T_n(f))\rangle$ by $(\sup_{\lambda}
1/h(\lambda))^2 \|D_n(S_nf)-T_n(f)\|_2^2$ is too rough (cf. also
Remark \ref{remondiffofgenmainthm}). Therefore, we look for a
new function, say $g$, with the properties that $1/(g h)$ behaves like
a constant for small $\lambda$ and $D_n(S_ng)$ is explicitly known. It
turns out that $g=f_{1/2+\alpha}$ is a good choice, where
$f_{1/2+\alpha}$ denotes the spectral density of a fractional Gaussian
noise process with Hurst index $1/2+\alpha\leq1/4$, cf. Lemma C.7 for
details. Furthermore, $r_{1/2+\alpha}$ denotes the corresponding
autocovariance function. From (\ref{eqfdmtlfXineqs}) and (\ref
{eqfdmtlhineqs}), we obtain
%
\begin{equation}\label{eqfhconstprop}
\sup_{\lambda\in[1/n, \uppi )} \frac1{f_{1/2+\alpha
}(\lambda) h(\lambda)}
\lesssim n^{2\beta+\varepsilon}.
\end{equation}
Define
\begin{eqnarray*}
\mathrm{I} &:=& \bigl\langle D_n^{-1} \bigl({hf_{1/2+\alpha}}\bigr)
\bigl(D_n (f_{1/2+\alpha} )-D_n (S_n
f_{1/2+\alpha} ) \bigr) \bigl(D_n (S_nf
)-T_n (f ) \bigr) \bigr\rangle,
\\
\mathrm{II} &:=& \bigl\langle\bigl(D_n (S_n f_{1/2+\alpha}
)-T_n (f_{1/2+\alpha} ) \bigr) \bigl(D_n
(S_nf )-T_n (f ) \bigr) \bigr\rangle,
\\
\mathrm{III} &:=& \bigl\langle T_n (f_{1/2+\alpha} ) \bigl(D_n
(S_nf )-T_n (f ) \bigr) \bigr\rangle
\end{eqnarray*}
and note that by Lemma C.1(iii),
%
\begin{eqnarray}\label{eqtoshowforalphasmall}
\tfrac14 \bigl\langle D_n^{-1} ({h}) \bigl(D_n
(S_nf )-T_n (f ) \bigr) \bigr\rangle&\leq& \mathrm{I} +
\bigl\|D_n^{-1} ({hf_{1/2+\alpha}}) \bigr\|_\infty^2
(\mathrm{II}+\mathrm{III}) \nonumber\\[-8pt]\\[-8pt]
&\lesssim& \mathrm{I} + n^{4\beta+2\varepsilon}
(\mathrm{II}+\mathrm{III}).\nonumber
\end{eqnarray}
In the following we will bound the terms $\mathrm{I}, \mathrm{II}$ and
$\mathrm{III}$, separately.

\mbox{}\hphantom{II}(I): By (ii) and (iv) of Lemma C.7, we know $(r_{1/2+\alpha
}(k))_k \in\GM$ and that $f_{1/2+\alpha}$ behaves like a multiple of
$\lambda^{-2\alpha}$ for $\lambda\downarrow0$. Using Lemma C.3(i)
and (\ref{eqfhconstprop}),
\begin{eqnarray*}
\sup_{\lambda\in[1/n, \uppi )} \frac{ |f_{1/2+\alpha}(\lambda) -S_n
f_{1/2+\alpha}(\lambda
) |}{ h(\lambda)f_{1/2+\alpha}(\lambda)} &\lesssim& n^{2\beta+1+\varepsilon}
\Biggl(\bigl|r_{1/2+\alpha}(n)\bigr|+\sum_{k=n+1}^\infty
\frac{|r_{1/2+\alpha}(k)|}{k} \Biggr) \\
&\lesssim& n^{2\beta+2\alpha
+\varepsilon}.
\end{eqnarray*}
Thus, with (\ref{eqdefdn}) and Lemma C.1(iii),
\begin{eqnarray*}
\mathrm{I}&\leq&\bigl\|D_n^{-1} ({hf_{1/2+\alpha}})
\bigl(D_n (f_{1/2+\alpha} )-D_n (S_n
f_{1/2+\alpha} ) \bigr) \bigr\| _\infty^2 \bigl\langle
D_n (S_nf )-T_n (f ) \bigr\rangle\\
&\lesssim&
n^{4\varepsilon
+4\beta}.
\end{eqnarray*}

\mbox{}\hphantom{I}(II): By Lemma C.7(i), and the boundedness of the sequence
$(r_{1/2+\alpha}(k))_k$, we see that there exists a constant
$C_{\alpha}$ such that for all $k\in\mathbb{N}$, $|r_{1/2+\alpha
}(k)|\leq C_\alpha k^{2\alpha-1}$. Using Lemma C.2 and (\ref
{eqfdmtlrXineqs}),
\begin{eqnarray*}
e_{i,j}:\!&=& \bigl| \bigl[ \bigl(D_n (S_n
f_{1/2+\alpha} )-T_n (S_n f_{1/2+\alpha} ) \bigr)
\bigl(D_n (S_nf )-T_n (f ) \bigr)
\bigr]_{i,j} \bigr|
\\
&\leq& C_\alpha C_r\sum_{k=1}^n
\bigl(|i+k|^{2\alpha
-1}+|2n+2-k-i|^{2\alpha-1} \bigr)\\
&&\hspace*{36pt}{}\times
\bigl(|j+k|^{-2\alpha-1+\delta'}+|2n+2-j-k|^{-2\alpha-1+\delta
'} \bigr).
\end{eqnarray*}
Define $F(i,j):= \sum_{k=1}^n |i+k|^{2\alpha-1} |j+k|^{-2\alpha-1}$.
It is well known that if $(a_k)_k$ and $(b_k)_k$ are nonnegative
sequences which are monotone increasing and decreasing, respectively,
then $\sum_{k=1}^n a_kb_k\leq\sum_{k=1}^n a_kb_{n+1-k}$. Thus,
\[
e_{i,j}\leq C_\alpha C_r(2n)^{\delta'}
\bigl(F(i,j)+F(i,n+1-j)+F(n+1-i,j)+F(n+1-i,n+1-j) \bigr).
\]
From the monotonicity of $x\mapsto x^{-2\alpha-1}$ and $x\mapsto
x^{2\alpha-1}$ for $x>0$,
\[
F(i,j) \leq j^{-2\alpha-1} \sum_{k=1}^n
|i+k|^{2\alpha-1} \leq j^{-2\alpha-1}\int_i^\infty
x^{2\alpha-1} \mrmd x = \frac
{1}{2|\alpha|} j^{-2\alpha-1}i^{2\alpha}.
\]
This allows to bound $e_{i,j}e_{j,i}$ by a multiple of $n^ {2\varepsilon}
(\min(i,n+1-i)\min(j,n+1-j))^{-1}$. Hence, $\mathrm{II} \leq\sum_{i,j=1}^n
e_{i,j}e_{j,i}\lesssim n^{2\varepsilon}\log^2n$.

(III): First let us introduce the projection $\Pi=\Pi_n$
defined for an $n\times n$ matrix $A=(a_{i,j})_{i,j=1,\ldots,n}$ by
$(\Pi A)_{i,j}:=a_{i,j}$ if $i+j\leq n+1$ and zero otherwise. Further
let $E$ denote the $n\times n$ matrix $(E)_{i,j}:=1$ if $i+j=n+1$ and
zero otherwise. In particular $E^2= \id_n$. Note that by Lemma C.1(iv),
\begin{eqnarray*}
&& \bigl\langle T_n (f_{1/2+\alpha} ) \bigl(D_n
(S_nf )-T_n (f ) \bigr) \bigr\rangle
\\
&&\quad\leq2 \bigl\langle T_n (f_{1/2+\alpha} ) \Pi
\bigl(D_n (S_nf )-T_n (f ) \bigr) \bigr\rangle
+ 2 \bigl\langle T_n (f_{1/2+\alpha} ) (\id_n-\Pi)
\bigl(D_n (S_nf )-T_n (f ) \bigr) \bigr
\rangle,
\end{eqnarray*}
where by a slight abuse of language $\id_n$ denotes here the identity
operator on the space of $n\times n$ matrices. To bound the first term,
we decompose
\begin{eqnarray*}
\bigl(T_n (f_{1/2+\alpha} ) \Pi\bigl(D_n
(S_nf )-T_n (f ) \bigr) \bigr)_{i,j} &=& \sum
_{k=1}^{2i-1} r_{1/2+\alpha}(i-k) [
\gamma_{k+j-1}-\gamma_{i+j-1} ]
\\
&&{} +\gamma_{i+j-1} \sum_{k=1}^{2i-1}
r_{1/2+\alpha}(i-k)
\\
&&{} +\sum_{k=2i}^{n+1-j} r_{1/2+\alpha}(i-k)
\gamma_{k+j-1}
\\
&=:&A_1(i,j)+A_2(i,j)+A_3(i,j)
\end{eqnarray*}
with the convention $\sum_{k=2i}^{n+1-j}=- \sum_{k=n+1-j}^{2i}$ if
$2i>n+1-j$. By assumption $|\gamma_{q+p}-\gamma_p|\leq\sum
_{v=p}^{q+p-1}|\gamma_{v+1}-\gamma_v|\lesssim n^\varepsilon qp^{-2\alpha
-1}$. Hence, uniformly in $i,j$,
%
\begin{equation}\label{eqA1bd}
\bigl|A_1(i,j)\bigr|\lesssim n^\varepsilon\sum
_{k=1}^{2i-1} \min\bigl(|i-k|^{2\alpha
},1\bigr)
j^{-2\alpha-2} \lesssim n^\varepsilon i^{2\alpha+1}j^{-2\alpha-2}.
\end{equation}
If $2j<i$, we can split the sum $\sum_{k=1}^{2i-1}= \sum_{k=1}^{\lfloor
i/2\rfloor}+\sum_{k=\lfloor i/2\rfloor+1}^{2i-1}$.
Then, the first part of $|A_1(i,j)|$ can be bounded also by
\[
\Biggl|\sum_{k=1}^{\lfloor i/2\rfloor} r_{1/2+\alpha}(i-k) [
\gamma_{k+j-1}-\gamma_{i+j-1} ] \Biggr| \lesssim n^\varepsilon
i^{2\alpha-1}\sum_{k=1}^{\lfloor i/2\rfloor}
(k+j-1)^{-2\alpha-1} \lesssim n^\varepsilon i^{-1}
\]
and the second part is by the same arguments as in (\ref{eqA1bd})
($j$ can now be replaced by $i+j$) of the order $n^\varepsilon i^{2\alpha
+1}(i+j-1)^{-2\alpha-2} \leq n^\varepsilon i^{-1}$. Together, this shows that
\[
\bigl|A_1(i,j)\bigr|\lesssim\cases{ n^\varepsilon i^{2\alpha+1}j^{-2\alpha-2},
&\quad if $2j\geq i$,
\cr
n^\varepsilon i^{-1}, &\quad if $2j<i$.}
\]
With Lemma C.7(i) and telescoping, $\sum_{k=1}^{2i-1} r_{1/2+\alpha
}(i-k) = i^{2\alpha+1}-(i-1)^{2\alpha+1}$ and therefore,
\[
\bigl|A_2(i,j)\bigr|\lesssim n^\varepsilon(i+j-1)^{-2\alpha-1}i^{2\alpha}.
\]
The last term of the expansion can be simply bounded by
\[
\bigl|A_3(i,j)\bigr|\lesssim n^\varepsilon(i+j-1)^{-2\alpha-1} \sum
_{k=2i}^{n} |i-k|^{2\alpha-1} \lesssim
n^\varepsilon(i+j-1)^{-2\alpha-1} i^{2\alpha}.
\]
In particular, the bounds for $A_2(i,j)$ and $A_3(i,j)$ are uniformly
in $i,j$ as well. Hence, by elementary computations
\begin{eqnarray*}
&& \bigl| \bigl\langle T_n (f_{1/2+\alpha} ) \bigl(D_n
(S_nf )-T_n (f ) \bigr) \bigr\rangle\bigr|
\\
&&\quad\leq\sum_{i,j} \bigl| A_1(i,j)+A_2(i,j)+A_3(i,j)
\bigr| \bigl| A_1(j,i)+A_2(j,i)+A_3(j,i) \bigr| \lesssim
n^{2\varepsilon} \log^2n.
\end{eqnarray*}
Finally, note that $E T_n (f_{1/2+\alpha} ) E = T_n
(f_{1/2+\alpha} )$
and $E ((\id_n-\Pi) [D_n (S_nf )-T_n (f
)] )E$ is a matrix
with entries $-\gamma_{i+j}$ for $i+j\leq n$ and zero otherwise.
Therefore, we have by rewriting
\[
\bigl\langle T_n (f_{1/2+\alpha} ) (\id_n-\Pi)
\bigl(D_n (S_nf )-T_n (f ) \bigr) \bigr\rangle
= \bigl\langle T_n (f_{1/2+\alpha} ) E \bigl[ (\id_n-
\Pi) \bigl(D_n (S_nf )-T_n (f ) \bigr)
\bigr]E \bigr\rangle
\]
the same structure as above (up to an index shift by one) and all
arguments apply. This shows that $\mathrm{III} \lesssim n^{2\varepsilon} \log^2 n$.

The estimates in (I), (II), and (III) show that the r.h.s. of
(\ref{eqtoshowforalphasmall}) can be upper bounded by
$n^{4\varepsilon+4\beta}\log^2 n$, and hence assumption (iii) of
Theorem \ref{thmgeneralizedmaintheorem} follows by choosing
$\varepsilon$ sufficiently small.

Since we have verified the assumptions of Theorem \ref
{thmgeneralizedmaintheorem}, part \ref{version3} of Theorem \ref
{mainthm} follows.
\end{pf*}
\end{appendix}

\section*{Acknowledgments}

Till Sabel was supported by DFG/SNF Grant FOR 916. Johannes
Schmidt-Hieber was funded by DFG postdoctoral fellowship SCHM 2807/1-1.
We would like to thank Axel Munk for helpful discussions.

\begin{supplement}
\stitle{Supplement to ``Asymptotically efficient estimation of a scale
parameter in Gaussian time series and closed-form expressions for the
Fisher information''}
\slink[doi]{10.3150/12-BEJ505SUPP} 
\sdatatype{.pdf}
\sfilename{BEJ505\_supp.pdf}
\sdescription{In the supplement, we provide the proofs of Theorem
\ref{thmexplicit} and Corollary \ref{cortheoremholdsforfBM} along
with some technical propositions and lemmas denoted by B.1, B.2$,\ldots,$C.1,
C.2$,\ldots$\,.}
\end{supplement}


\printhistory

\end{document}